\documentclass[12pt,a4paper]{amsart}

\usepackage{amssymb,amscd,amsthm,amsxtra,mathtools,amsmath}
\usepackage[margin=2cm, top=3cm]{geometry}
\usepackage[alphabetic]{amsrefs}
\usepackage{hyperref}
\usepackage{mathtools}
\usepackage{xcolor}
\usepackage{graphicx}
\usepackage{wrapfig}
\usepackage{float}
\usepackage{enumitem}
\usepackage{mathrsfs}
\usepackage{comment}

\numberwithin{equation}{section}
\allowdisplaybreaks[4]

\newtheorem{theorem}{Theorem}[section]
\newtheorem{proposition}[theorem]{Proposition}

\newtheorem{lemma}[theorem]{Lemma}
\newtheorem{question}[theorem]{Question}
\theoremstyle{definition}

\theoremstyle{remark}
\newtheorem{remark}[theorem]{Remark}

\newcommand{\dist}{\operatorname{dist}}
\newcommand{\R}{\mathbb{R}}
\newcommand{\N}{\mathbb{N}}
\newcommand{\eps}{\varepsilon}

\newcommand{\restr}[2]{{\left.\kern-\nulldelimiterspace#1\vphantom{\big|}\right|_{#2}}}

\newcommand{\per}{\operatorname{Per}}
\newcommand{\pers}{\operatorname{Per}_s}
\newcommand{\pv}{\operatorname{p.v.}}
\newcommand{\loc}{\operatorname{loc}}
\newcommand{\sep}{\operatorname{sep}}

\renewcommand{\leq}{\leqslant}
\renewcommand{\le}{\leqslant}
\renewcommand{\geq}{\geqslant}

\renewcommand{\epsilon}{\varepsilon}
\renewcommand{\O}[2]{\mathcal{O}^{(#1)} (#2)}

\newcommand{\pd}[2]{#1_{,#2}}         % first-order partial derivative,  e.g. \pd{u}{j}, \pd{\nu}{j}
\newcommand{\pdd}[3]{#1_{,#2#3}}      % second-order partial derivative, e.g. \pdd{u}{i}{j}
\newcommand{\hess}[3]{#1_{;#2#3}}     % covariant Hessian,               e.g. \hess{u}{i}{j}
\renewcommand{\pd}[2]{{\partial_{#2}#1}}
\renewcommand{\pdd}[3]{{\partial_{#2#3}#1}}
\renewcommand{\hess}[3]{{\nabla^2_{#2#3}#1}}

\mathtoolsset{showonlyrefs=true}
\begin{document}
\title[Nonlocal minimal surfaces with~$k$ connected components]{Nonlocal approximation of stable minimal surfaces with prescribed number of connected components}
\author[L. De Gaspari]{Lorenzo De Gaspari}
\address{L. D. G., 
Department of Mathematics and Statistics,
University of Western Australia,
35~Stirling Highway, WA 6009 Crawley, Australia.}

\email{lorenzo.degaspari@research.uwa.edu.au}

\author[S. Dipierro]{Serena Dipierro}
\address{S. D., 
Department of Mathematics and Statistics,
University of Western Australia,
35~Stirling Highway, WA 6009 Crawley, Australia.}

\email{serena.dipierro@uwa.edu.au}

\author[O. Savin]{Ovidiu Savin}
\address{O. S., 
Department of Mathematics,
Columbia University,
2990 Broadway, NY 10027 New York, US.}

\email{savin@math.columbia.edu}

\author[E. Valdinoci]{Enrico Valdinoci}
\address{E. V., 
Department of Mathematics and Statistics,
University of Western Australia,
35~Stirling Highway, WA 6009 Crawley, Australia.}

\email{enrico.valdinoci@uwa.edu.au}

\thanks{LDG is supported by a Scholarship for International Research Fees at the University of Western Australia. LDG and EV are supported by Australian Research Council FL190100081. SD is supported by Australian Research Council FT230100333. SD and EV are members of the Australian Mathematical Society (AustMS). OS is supported by the NSF Grant DMS-2349794.}



\begin{abstract}
    We show that every strictly stable classical minimal hypersurface in a Euclidean ball can be approximated, for~$s$ sufficiently close to~$1$, by a stable~$s$-minimal surface in a smaller ball. We also give an explicit lower bound on the size of the interval of admissible values of~$s$.
    
    As a consequence, we construct stable nonlocal minimal surfaces with an arbitrarily large number~$k$ of connected components in a bounded region, under the quantitative condition~$k \leq \frac{C}{\sqrt{1-s}}$, attaining an optimal result (at least in dimensions~$2$ and~$3$). 

    As part of the proof of our main result, we obtain an auxiliary result of independent interest, providing quantitative estimates for foliations of constant mean curvature embedded hypersurfaces near a stable minimal surface. The resulting foliation supplies nonlocal barriers thanks to known results of uniform convergence of the nonlocal mean curvature to the classical one when~$s \to 1$.
\end{abstract}
\maketitle
\section{Introduction}
A central research question in the area of partial differential equations and geometric analysis is to understand the connection between classical and nonlocal minimal surfaces. Developing a deeper understanding of the relationship between the classical and fractional settings is of fundamental importance, as it may allow results and methodologies to be transferred from one theory to the other, combining the extensive body of knowledge available for classical problems with the novel techniques and phenomena arising in the nonlocal setting.

For instance, classical minimal surfaces benefit from a well-developed regularity theory, whereas the regularity theory for nonlocal minimal surfaces still presents several fundamental open problems (see, e.g.,~\cite{MR3824212}). Conversely, nonlocal minimal surfaces exhibit distinctive rigidity phenomena that have no counterpart in the classical theory (see, e.g.,~\cite{MR3981295}).

A particularly interesting class is provided by stable surfaces, namely surfaces with vanishing first variation of the (non)local perimeter functional and nonnegative second variation. From the viewpoint of geometric complexity, it is especially intriguing to approximate a given classical surface by a large number of nonlocal sheets. Such an approximation is known in the case of hyperplanes (see~\cite{MR4116635}), while upper bounds on the number of sheets have been established in low dimensions (see~\cite{chan2025nonlocalapproximationminimalsurfaces}). %chktex 36

In this paper, we address this problem by constructing arbitrarily many stable nonlocal sheets in a prescribed neighborhood of a stable classical minimal surface. We further establish a sharp bound for the number of such sheets in terms of the fractional parameter. More specifically, we are interested in understanding stable nonlocal minimal surfaces as small perturbations of stable classical minimal surfaces in the asymptotic regime as~$s \to 1$. This leads to the following questions.

\begin{question}\label{QST::Perturbation}
    Suppose that~$\partial E$ is a strictly stable minimal hypersurface. Can we find~$s \in (0,1)$ `close' to~$1$ and a stable~$s$-minimal hypersurface~$\partial F$ `close' to~$\partial E$?
\end{question}

\begin{question}\label{QST::Sheeting}
    Let~$k \in\N\setminus\{0\} $. Can we find~$s \in (0,1)$ such that~$\frac{1}{\sqrt{1-s}}$ is comparable to~$k$, and a relative minimizer~$F$ of the~$s$-perimeter such that~$\partial F$ has~$k$ connected components in a bounded region?
\end{question}

Let us precisely introduce the mathematical setting in which we will provide an answer to these informal questions.

\subsection{Classical minimal surfaces}
Following the variational approach of De Giorgi~\cites{MR62214,MR74499} (see also~\cite{MR775682}), for an open set~$\Omega \subset \R^n$ with~$n \geq 2$, we identify a minimal surface with the boundary of a measurable set~$E$ that is a critical point of the perimeter functional
\begin{equation}\label{EQ::PerimeterDefinition}
    \per(E, \Omega) \coloneq \sup\left\{\int_{E} \operatorname{div}\phi\,dx\colon \phi \in C^1_c(\Omega; \R^n),\,\|\phi\|_{L^\infty(\Omega)}\leq 1\right\}.
\end{equation}
As customary, see e.g.~\cite{MR3516886}*{Appendix~A}, the sets we consider are always assumed to be open in the measure theoretic sense. In this setting, the minimality condition for~$E$ is defined with respect to compactly supported deformations. Namely, for every~$X \in C^\infty_c(\Omega; \R^n)$, we denote by~$\Phi^X\colon \Omega \times \R \to \Omega$ the associated vector flow, i.e., the unique solution of the Cauchy problem
\begin{equation}
    \left\{
        \begin{aligned}
            &\partial_t \Phi^X(x,t) = X(\Phi^X(x,t))&\text{for every }x \in \Omega \text{ and } t > 0,\\
            &\Phi^X(x, 0) = x&\text{for every }x \in \Omega.
        \end{aligned}
    \right.
\end{equation}
Then,~$\partial E$ is minimal in~$\Omega$ if
\begin{equation}\label{EQ::MinimalSurfaceFlowDefinition}
    \restr{\frac{d}{dt}}{t=0}\per(\Phi^X(E, t), \Omega) = 0
\end{equation}
for every~$X \in C^\infty_c(\Omega; \R^n)$. We stress that, under this definition, a minimal surface is not necessarily a minimizer of the perimeter functional. Whenever~$\partial E \cap \Omega$ is a~$C^2$ submanifold of~$\R^n$, minimality is equivalent (see e.g.~\cites{MR775682,MR2780140}) to the pointwise satisfaction of the Euler-Lagrange equation
\begin{equation}\label{EQ::ClassicalEulerLagrange}
    H[\partial E](x) = 0,
\end{equation}
where~$H[\partial E](x)$ denotes the mean curvature, defined as the sum of the principal curvatures of~$\partial E$ at~$x$.
Here and throughout the paper, the principal curvatures are the eigenvalues of the shape operator~$d\nu$, where~$\nu$ is the outer unit normal to~$E$ (or, for an orientable hypersurface~$\Sigma$, the chosen unit normal). With this convention, the boundary of a ball has positive mean curvature. If a~$C^2$ orientable hypersurface~$\Sigma \subset \R^n$ has vanishing mean curvature in its relative interior, we may also refer to~$\Sigma$ as a minimal hypersurface\footnote{Note that the two notions coincide in the following sense. Whenever~$\partial E$ is a~$C^1$ submanifold, we have that~$\per(E, \Omega) = \mathcal{H}^{n-1}(\partial E \cap \Omega)$. One can prove that
\begin{equation}
    \restr{\frac{d}{dt}}{t=0}\mathcal{H}^{n-1}(\Phi^X(\Sigma, t)) = 0
\end{equation}
whenever~$\Sigma$ is minimal and~$X$ is a smooth vector field that vanishes on~$\partial \Sigma$.
}.

An important property in the study of minimal surfaces is that of stability. Informally, a stable minimal surface is a minimizer of the perimeter `up to second order'. In our framework we say that~$\partial E$ is a stable minimal surface if~\eqref{EQ::MinimalSurfaceFlowDefinition} is satisfied and
\begin{equation}\label{EQ::StabilityFlowDefinition}
    \restr{\frac{d^2}{dt^2}}{t=0}\per(\Phi^X(E, t), \Omega) \geq 0.
\end{equation}
The stability of a minimal hypersurface~$\Sigma$ is also characterized by the spectral stability of the Jacobi\footnote{In the more general case of minimal submanifolds of a Riemannian manifold, the Jacobi operator contains an additional term given by the ambient Ricci curvature. Clearly, this term vanishes whenever the ambient space is Euclidean.} operator
\begin{equation}\label{EQ::JacobiOperatorDefinition}
\mathcal{L}_\Sigma \coloneq -\Delta_{\Sigma} - |A|^2
\end{equation}
(see e.g.,~\cite{MR562550}*{Section~3} and the references therein, or~\cite{MR2780140}*{Proposition~1.39}). Namely, denoting by~$\Delta_{\Sigma}$ the Laplace-Beltrami operator on~$\Sigma$ and by~$A$ the second fundamental form, we have that~$\Sigma$ is a stable critical point of the perimeter if and only if
\begin{equation}\label{EQ::StabilitySpectralDefinition}
    \lambda_1(\Sigma) \coloneq \inf_{\substack{\phi \in C^{\infty}_c(\Sigma)\\\phi \not\equiv 0}}\frac{\displaystyle\int_{\Sigma} \phi \mathcal{L}_\Sigma \phi\,d\mathcal{H}^{n-1}}{\displaystyle\int_{\Sigma} \phi^2\,d\mathcal{H}^{n-1}} \geq 0.
\end{equation}
If the inequality holds strictly in~\eqref{EQ::StabilityFlowDefinition} (respectively,~\eqref{EQ::StabilitySpectralDefinition}), we say that~$\partial E$ (respectively,~$\Sigma$) is a strictly stable minimal surface. We remark that a strictly stable critical point is a local minimizer, although the converse is not true in general.

\subsection{Nonlocal minimal surfaces}
Nonlocal minimal surfaces were introduced in~\cite{MR2675483} as the boundaries of sets that are critical points of the fractional perimeter, defined by taking~$s \in (0,1)$ and letting
\begin{equation}\label{EQ::sPerimeterDefinition}
    \pers(E, \Omega) \coloneq L(E \cap \Omega, E^c \cap \Omega) + L(E \cap \Omega, E^c \cap \Omega^c) + L(E \cap \Omega^c, E^c \cap \Omega),
\end{equation}
where we use the superscript letter~$c$ to denote the complement of a set in~$\R^n$, and~$L(\cdot,\cdot)$ is the nonlocal interaction
\begin{equation}
    L(A,B) \coloneq \int_{\R^n \times \R^n} \frac{\chi_A(x) \chi_B(y)}{|x-y|^{n+s}}\,dxdy.
\end{equation}

Although the definitions~\eqref{EQ::PerimeterDefinition} and~\eqref{EQ::sPerimeterDefinition} may look unrelated at face value, as~$s \to 1$, the~$s$-perimeter converges to the classical perimeter under different assumptions and notions of convergence, see~\cites{MR2765717,MR2782803,MR3107529}. In fact, both classical and nonlocal minimal surfaces are related to phase transition models in physics, where taking into account long range effects naturally leads to nonlocality, see also~\cite{MR4581189} and the references therein.

As in the classical case, we say that an~$s$-minimal surface is a critical point of the~$s$-perimeter. That is,~$\partial E$ is~$s$-minimal in~$\Omega$ if, for every~$X\in C^\infty_c(\Omega; \R^n)$, we have that
\begin{equation}
    \restr{\frac{d}{dt}}{t=0}\pers(\Phi^X(E, t), \Omega) = 0,
\end{equation}
and stability is given by
\begin{equation}
    \restr{\frac{d^2}{dt^2}}{t=0}\pers(\Phi^X(E, t), \Omega) \geq 0.
\end{equation}
If the inequality is strict, we say that~$\partial E$ is strictly stable.

The Euler-Lagrange equation for a critical point~$E$ of the~$s$-perimeter, posed on~$\partial E$, reads
\begin{equation}\label{EQ::NonlocalEulerLagrange}
    H_s[\partial E](x) \coloneq \pv\int_{\R^n} \frac{ \chi_{E^c}(y) - \chi_{E}(y)}{|x-y|^{n+s}}\,dy = 0,
\end{equation}
where the symbol~$\pv$ indicates that the integral is taken in the principal value sense. By analogy with the classical case~\eqref{EQ::ClassicalEulerLagrange}, the quantity~$H_s$ defined by~\eqref{EQ::NonlocalEulerLagrange} is often referred to as the nonlocal mean curvature. Notice that, with this sign convention, the boundary of a ball has positive nonlocal mean curvature, consistently with the classical case.
See also~\cite{MR3230079} for additional remarks about analogies between these two notions of mean curvature.

\subsection{Statement of the main results}
Throughout this paper, we denote by~$B_r(x)$ the ball of radius~$r > 0$ centered at the point~$x \in \R^n$, and we may omit the dependence on~$x$ whenever~$x = 0$ by simply writing~$B_r$ in place of~$B_r(0)$.

\begin{figure}[ht]
\centering
\includegraphics[width=0.8\textwidth]{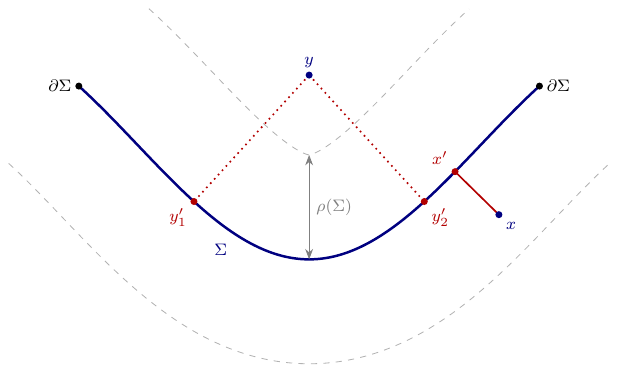}
\caption{Sketch of a hypersurface with positive reach. The point~$x$ lies within the reach and admits the unique projection~$x^\prime$. On the contrary, the point~$y$ has a distance from~$\Sigma$ greater than~$\rho(\Sigma)$, and it admits two different projections~$y^\prime_1$ and~$y^\prime_2$. }\label{FIG::Reach}
\end{figure}

For a compact, orientable, embedded~$C^2$ hypersurface~$\Sigma$, we denote by~$\rho(\Sigma)$ the reach of~$\Sigma$ in the sense of Federer~\cite{MR110078}. That is,~$\rho(\Sigma)$ is the largest positive value such that, for every~$x \in \R^n$ satisfying~$\operatorname{dist}(x,\Sigma) < \rho(\Sigma)$, there exists a unique projection~$x^\prime \in \Sigma$ attaining the optimal distance~$|x-x^\prime| = \operatorname{dist}(x,\Sigma)$. See also Figure~\ref{FIG::Reach} for an intuitive representation. In the case where~$\Sigma$ is the union of~$k \in \N\setminus\{0\}$ connected components~${\{\Sigma^j\}}_{j=1}^k$, we make an important distinction between two different aspects that bound~$\rho(\Sigma)$ from above. On one hand, it is clear that
\begin{equation}\label{EQ::rhoLocDefinition}
    \rho_{\loc}(\Sigma) \coloneq \min_{1\leq j \leq k} \rho(\Sigma^j) \geq \rho(\Sigma).
\end{equation}
On the other hand, we have that the separation
\begin{equation}
    \sep(\Sigma) \coloneq \frac{1}{2}\min_{1\leq i < j \leq k} \dist(\Sigma^i, \Sigma^j) \geq \rho(\Sigma),
\end{equation}
with the convention that 
if~$k=1$ then~$\sep(\Sigma)=+\infty$.
Accordingly, if~$\Sigma$ is a connected hypersurface, then~$\rho(\Sigma) = \rho_{\loc}(\Sigma)$ as expected.

The notation and terminology that we introduced up to now allow us to state the main result of the present paper.

\begin{theorem}\label{THM::Main}
    Let~$k \in \N\setminus\{0\}$,~$\alpha \in (0,1)$, and~$E \subset \R^n$. Suppose that~$\partial E \cap B_2$ is an orientable,~$C^{3,\alpha}$, strictly stable minimal hypersurface in~$B_2$. Suppose also that~$\partial E \cap B_{\frac{1}{2}}$ has~$k$ connected components. Assume that the~$C^{3,\alpha}$ norm of~$\partial E$ is bounded by~$M > 0$.
    
    Then, there exists a constant~$C > 0$, which depends only on~$n$, $\alpha$, $\lambda_1(\partial E \cap B_2)$, $M$, and~$\rho_{\loc}(\partial E \cap B_2)$, such that, for every
    \begin{equation}\label{EQ::s0Constraint}
        s \in \left[1-C\min\{1,\sep^2(\partial E \cap B_2)\}, 1\right), %chktex 9
    \end{equation}
    there exists~$F \subset \R^n$ such that~$\partial F$ is stable and~$s$-minimal in~$B_1$, and~$\partial F \cap B_1$ has at least~$k$ connected components. %chktex 9
\end{theorem}

Our proof of Theorem~\ref{THM::Main} relies on a Perron-type method. The fundamental ingredient is a quantitative foliation result that enables the construction of the necessary barriers. Such result is stated as follows.

\begin{theorem}\label{THM::Foliation}
    Let~$\alpha \in (0,1)$ and let~$\Sigma \subset B_2$ be a compact, orientable,~$C^{3,\alpha}$ hypersurface embedded in~$\R^n$, with boundary~$\partial \Sigma = \Sigma \cap \partial B_2 \neq \varnothing$. Suppose that~$\Sigma$ is a strictly stable minimal surface in~$B_2$. Assume that the~$C^{3,\alpha}$ norm of~$\Sigma$ is bounded by~$M > 0$.

    Then, there exist~$\eps_0 > 0$ and a one-parameter family of embedded hypersurfaces~${\{\Sigma_t\}}_{t \in (-\eps_0, \eps_0)}$, whose elements foliate a neighborhood of~$\Sigma \smallsetminus \partial \Sigma$ such that~$\Sigma_0 = \Sigma$ and, for every~$t \in (-\eps_0, \eps_0)$, we have that
    \begin{equation}\label{EQ::FoliationCMCCondition}
        H[\Sigma_t] \equiv t.
    \end{equation}

    Moreover, there exists a constant~$c > 0$, which depends only on~$n$,~$\alpha$,~$\lambda_1(\Sigma)$,~$M$, and~$\rho_{\loc}(\Sigma)$, such that the 
    largest possible value of~$\eps_0$ for which the above
    statement holds true satisfies
    \begin{equation}\label{EQ::eps0OptimalCondition}
        \eps_0 \geq c \min\{1,\sep(\Sigma)\}.
    \end{equation}
\end{theorem}

\begin{figure}[ht]
\centering
\includegraphics[width=0.8\textwidth]{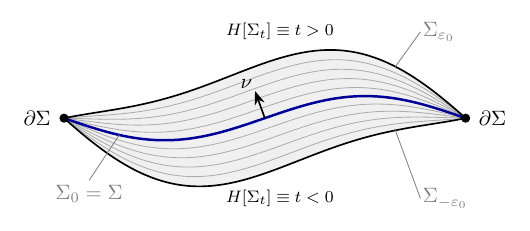}
\caption{Sketch illustrating the foliation around a connected hypersurface~$\Sigma$, as provided by Theorem~\ref{THM::Foliation}.}\label{FIG::FoliationConnected}
\end{figure}

\begin{figure}[ht]
\centering
\includegraphics[width=0.8\textwidth]{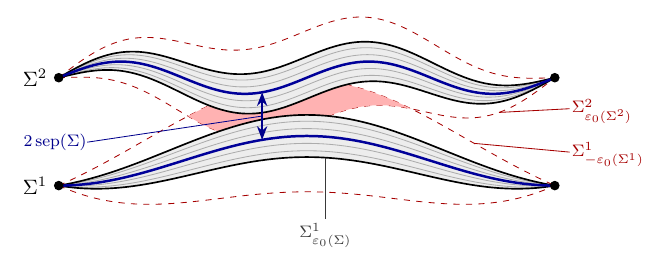}
\caption{Sketch illustrating the foliation around a hypersurface~$\Sigma$ with~$k = 2$ connected components, as provided by Theorem~\ref{THM::Foliation}.}\label{FIG::Foliation}
\end{figure}

See Figures~\ref{FIG::FoliationConnected} and~\ref{FIG::Foliation} for a visualization of the foliation provided by Theorem~\ref{THM::Foliation}. All the figures throughout this paper are schematic, and are only intended as visual aids to illustrate the geometric constructions.

\subsection{Comparison with previous results}
Theorem~\ref{THM::Main} provides an affirmative answer for Questions~\ref{QST::Perturbation} and~\ref{QST::Sheeting}. Note that when~$k$ is large,~$\sep(\partial E \cap B_2)$ is bounded from above, up to a multiplicative constant, by~$\frac{1}{k}$ (the bound is attained by considering~$k$ evenly spaced sheets), since we are fitting many connected components in the bounded region~$B_{\frac{1}{2}}$. This ensures that the lower bound for~$s$ in~\eqref{EQ::s0Constraint} is equivalent to~$k \leq \frac{C}{\sqrt{1-s}}$. In particular,~\cite{chan2025nonlocalapproximationminimalsurfaces}*{Theorem~1.1} yields the optimality of our result in dimensions~$2$ and~$3$. We believe that our result may be optimal also in higher dimension.

Our result may also be compared to~\cite{MR3798717}*{Theorem~2}, where the authors find a global nonlocal minimal surface with multiple `sheets'. Our result provides stronger quantitative estimates and hints towards a stronger genericity property of multi-sheet nonlocal minimal surfaces, at the price of working only at a local scale.

The proof of Theorem~\ref{THM::Foliation} is based on a strategy similar to the one employed by the authors of~\cite{MR2481045} to prove the existence of constant mean curvature submanifolds by perturbing geodesic spheres in a Riemannian manifold. Other similar approaches have been successfully employed in the literature to obtain comparable results, see for instance~\cite{MR1084717}. To the best of our knowledge, the main idea behind this structure of proof was introduced in~\cite{MR1101226}. In this work, the author proposes an abstract framework for the search of critical points of a parameter-dependent functional, applying this theory to the case of the mean curvature of an immersed submanifold in dependence of the chosen metric for the ambient manifold.

\subsection{Organization of the paper}

The rest of this paper is organized as follows. In Section~\ref{SEC::cmcFoliation} we discuss the effect of perturbations on the mean curvature of a regular hypersurface and we apply such results to prove Theorem~\ref{THM::Foliation} by employing a fixed point approach. Section~\ref{SEC::sMinimalSurface} contains the barrier construction and the Perron-type argument for the proof of Theorem~\ref{THM::Main}. In the Appendix, we show an explicit example in the flat case for Theorem~\ref{THM::Foliation}. 

\section{Construction of a CMC foliation and Proof of Theorem~\ref{THM::Foliation}}\label{SEC::cmcFoliation}
In this section we provide our proof of Theorem~\ref{THM::Foliation}, showing that the neighborhood of a smooth, classical minimal hypersurface is foliated by a family of constant mean curvature hypersurfaces. This construction provides the geometric backbone for the proof of Theorem~\ref{THM::Main} contained in Section~\ref{SEC::sMinimalSurface}.

Our proof strategy combines quantitative curvature estimates for~$C^{2,\alpha}$ perturbations of the hypersurface~$\Sigma$ together with a fixed point argument based on the contraction principle.

\subsection{Curvature estimates for perturbed hypersurfaces}
Here we collect some estimates regarding certain geometric quantities and their dependence on~$C^{2,\alpha}$ perturbations of a hypersurface. We provide the detailed computations for the convenience of the reader, but also because great care is necessary in order to obtain the precise relation between constants and some of the problem parameters.

\begin{figure}[ht]
\centering
\includegraphics[width=0.8\textwidth]{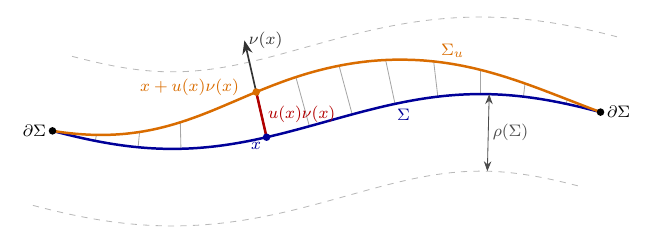}
\caption{Sketch of the normal geodesic graph~$\Sigma_u$ of a function~$u$ over a hypersurface~$\Sigma$.}\label{FIG::normalGraph}
\end{figure}

Consider a compact,~$C^{3,\alpha}$, orientable hypersurface~$\Sigma$. Denoting by~$\nu$ the unit normal vector field induced by the orientation, for a function~$u \in C^{2,\alpha}(\Sigma)$, we have that the normal geodesic graph of~$u$ over~$\Sigma$, given by
\begin{equation}\label{EQ::SigmaUDefinition}
    \Sigma_u \coloneq \{x + u(x)\nu(x) \colon x \in \Sigma\},
\end{equation}
is a~$C^{2,\alpha}$ hypersurface (see also Figure~\ref{FIG::normalGraph}) provided that the maximal absolute value attained by~$u(x)$ is strictly smaller than~$\rho(\Sigma)$.

In the case where~$\Sigma$ has nonempty boundary, we are interested in considering only variations that vanish on~$\partial \Sigma$, therefore we introduce the class of admissible functions
\begin{equation}\label{EQ::FSigmaDefinition}
    \mathcal{F}(\Sigma) \coloneq \left\{u \in C^{2,\alpha}(\Sigma)\colon u=0\text{ on }\partial\Sigma\text{ and }\|u\|_{C^{2,\alpha}(\Sigma)} < \frac{1}{2}\rho(\Sigma)\right\},
\end{equation}
which is an open subset of the Banach space of~$C^{2,\alpha}$ functions on~$\Sigma$ that vanish on the boundary, which we shall denote by~$C^{2,\alpha}_0(\Sigma)$.

We stress that, if~$u \in \mathcal{F}(\Sigma)$, then
\[\sup_{x \in \Sigma} |u(x)|\leq \|u\|_{C^{2,\alpha}(\Sigma)}< \frac{1}{2}\rho(\Sigma),\]
and therefore the map~$\Phi_u \colon \Sigma \to \Sigma_u$, defined as
\begin{equation}
    \Phi_u \colon x \longmapsto x + u(x)\nu(x),
\end{equation}
is a~$C^{2,\alpha}$ diffeomorphism. We analyze the dependence on~$u$ of the mean curvature of~$\Sigma_u$ by leveraging the pullback through~$\Phi_u$, which allows us to carry out the analysis on~$\Sigma$ itself.

We now fix some notation to improve the readability of formulas. We always implicitly assume the Einstein convention for sums over repeated indices. We denote by~$\pd{u}{j}$ and~$\pdd{u}{i}{j}$ the components with respect to local coordinates of the partial derivatives of~$u$ of first and second order respectively. We may also write~$\hess{u}{i}{j} \coloneq \pdd{u}{i}{j} - \Gamma_{ij}^k \pd{u}{k}$ for the covariant Hessian of~$u$, where~$\Gamma_{ij}^k$ are the Christoffel symbols of~$\Sigma$. For a given Riemannian metric~$g$ of components~$g_{ij}$, we shall denote by~$g^{ij}$ the components of the inverse of~$g_{ij}$, in the sense that~$g_{ij}g^{jk} \equiv \delta_{i}^k$, with~$\delta_i^k$ representing the Kronecker delta symbol.

In the wake of~\cite{MR2481045}, we introduce a compact notation for Taylor expansions of quantities that depend on a function. For any~$k \in \N\setminus\{0\}$, 
we let~$X$ be a nonlinear operator
in the function~$u \in \mathcal{F}(\Sigma)$ and its covariant derivatives up to the second order and we say that~$X=\O{k}{u}$ if~$X(0) = 0$ and there exists a constant~$C = C(n,\alpha,\|\Sigma\|_{C^{3,\alpha}}, \rho_{\loc}) > 0$ such that
\begin{equation}\label{EQ::OkNotationDefinition}
    \|X(u_1) - X(u_2)\|_{C^{0,\alpha}(\Sigma)} \leq C{\left(\|u_1\|_{C^{2,\alpha}(\Sigma)} + \|u_2\|_{C^{2,\alpha}(\Sigma)}\right)}^{k-1}\|u_1 - u_2\|_{C^{2,\alpha}(\Sigma)}
\end{equation}
We remark that, due to~\eqref{EQ::rhoLocDefinition} and~\eqref{EQ::FSigmaDefinition}, if~$X=\O{k+1}{u}$, then~$X = \O{k}{u}$. Also,\[\O{k}{u} + \O{m}{u} = \O{\min\{k,m\}}{u}\] and\[\O{k}{u}\O{m}{u} = \O{k+m}{u}.\]

We give an explicit expansion of the Riemannian metric of~$\Sigma_u$, pulled back to~$\Sigma$. We want to emphasize that, unlike the Taylor approximations given by the forthcoming Lemma~\ref{LEM::hTildeExpansion} and Proposition~\ref{PROP::CurvatureExpansion}, the perturbed metric that we compute does not contain an error term of type~$\O{k}{u}$.

\begin{lemma}\label{LEM::gTildeExpansion}
    Let~$\Sigma$ be a~$C^{3,\alpha}$, compact, orientable hypersurface embedded in~$\R^n$ and~$u \in \mathcal{F}(\Sigma)$. Let~$x \coloneq (x^1,\dots,x^{n-1})$ be local coordinates on a relatively open set~$U \subset \Sigma$, and~$\tilde{x} \coloneq x \circ \Phi_u^{-1}$ be the induced coordinates on~$\Phi_u(U)$.
    
    Denote by~$g_{ij}$ and~$\tilde{g}_{ij}$ the components of the Riemannian metrics on~$U$ and~$\Phi_u(U)$ induced by the ambient Euclidean metric, relative to the coordinates~$x$ and~$\tilde{x}$ respectively.
    
    Then, at every point~$p \in U$ and for every~$1\leq i,j\leq n-1$, the following expansion holds
    \begin{equation}\label{EQ::gTildeExpansion}
        {(\Phi_u^\ast \tilde{g})}_{ij} = g_{ij} + 2 u h_{ij} + u^2 h_i^k h_{jk} + \pd{u}{i} \pd{u}{j},
    \end{equation}
    where~$h_{i}^k$ denotes the shape operator of~$\Sigma$, and~$h_{ij} \coloneq g_{ik}h_j^k$ is the second fundamental form.
\end{lemma}
\begin{proof}
    Let~$p \in U$ and~$\partial_j \coloneq \restr{\frac{\partial}{\partial x_j}}{p}$ be the local frame on~$U$. We compute the differential of the map~$\Phi_u$, evaluated at each basis vector~$\partial_j$ as
    \begin{equation}\label{EQ::dPhifExplicit}
        d\Phi_u[\partial_j] = dx[\partial_j] + u d\nu[\partial_j]+ d u[\partial_j] \nu = \partial_j + u h_j^l \partial_l + \pd{u}{j}\nu,
    \end{equation}
    where we have also used the Weingarten equation~$\pd{\nu}{j} = h_j^l \partial_l$.

    Since~$\tilde{g}$ is induced by the Euclidean metric, the pullback metric is given by
    \begin{equation}\label{EQ::PullbackMetricEuclidean}
            {(\Phi_u^\ast \tilde{g})}_{ij} = \tilde{g}(d\Phi_u[\partial_i], d\Phi_u[\partial_j]) = \langle d\Phi_u[\partial_i], d\Phi_u[\partial_j] \rangle,
    \end{equation}
    where~$\langle\cdot,\cdot\rangle$ is the Euclidean scalar product. Substituting~\eqref{EQ::dPhifExplicit} into~\eqref{EQ::PullbackMetricEuclidean}, we find that
    \begin{equation}\label{EQ::PullbackMetricFullComputation}
        \begin{split}
            {(\Phi_u^\ast \tilde{g})}_{ij} &= \langle \partial_i, \partial_j \rangle + u h_i^k \langle \partial_k, \partial_j \rangle + \pd{u}{i} \langle \nu, \partial_j \rangle + u h_j^l \langle \partial_i, \partial_l \rangle + u^2 h_i^k h_j^l \langle \partial_k, \partial_l \rangle + \pd{u}{i}u h_j^l \langle \nu, \partial_l \rangle\\
            &\qquad + \pd{u}{j}\langle \partial_i, \nu \rangle + uh_i^k \pd{u}{j}\langle \partial_k, \nu \rangle + \pd{u}{i}\pd{u}{j}\langle \nu, \nu\rangle.
        \end{split}
    \end{equation}

    We recall that~$\langle \partial_i, \partial_j \rangle = g(\partial_i, \partial_j) = g_{ij}$ because~$g$ is induced by the Euclidean metric. Moreover, since~$\nu \perp T_p\Sigma$, we have~$\langle \nu, \partial_j \rangle = 0$. By assumption,~$|\nu|^2 = \langle \nu, \nu \rangle = 1$. All in all, these considerations allow the expression in~\eqref{EQ::PullbackMetricFullComputation} to be simplified, yielding
    \begin{equation}\label{EQ::PullbackMetricSimplified}
        {(\Phi_u^\ast \tilde{g})}_{ij} = g_{ij} + u h_{i}^k g_{kj} + u h_{j}^l g_{il} + u^2 h_i^k h_{jk} + \pd{u}{i} \pd{u}{j}.
    \end{equation}
    The desired result follows from~\eqref{EQ::PullbackMetricSimplified} upon noticing that~$h_{i}^k g_{kj} = h_{ij} = h_{ji} = h_{j}^l g_{il}$ due to the symmetry of the second fundamental form.
\end{proof}

We also find an expansion up to second order for the pulled back second fundamental form.

\begin{lemma}\label{LEM::hTildeExpansion}
    Let~$\Sigma$ be a~$C^{3,\alpha}$, compact, orientable hypersurface embedded in~$\R^n$ and~$u \in \mathcal{F}(\Sigma)$. Let~$x \coloneq (x^1,\dots,x^{n-1})$ be local coordinates on a relatively open set~$U \subset \Sigma$, and~$\tilde{x} \coloneq x \circ \Phi_u^{-1}$ be the induced coordinates on~$\Phi_u(U)$.
    
    Denote by~$h_{ij}$ and~$\tilde{h}_{ij}$ the components of the second fundamental form on~$U$ and~$\Phi_u(U)$, relative to the coordinates~$x$ and~$\tilde{x}$ respectively.

    Then, at every point~$p \in U$ and for every~$1\leq i,j\leq n-1$, the following expansion holds
    \begin{equation}\label{EQ::hTildeExpansion}
        {(\Phi_u^\ast \tilde{h})}_{ij} = h_{ij} - \hess{u}{i}{j} + u h_i^k h_{kj} + \O{2}{u}.
    \end{equation}
\end{lemma}
\begin{proof}
    Let~$p \in U$. We denote by~$\tilde{\nu}$ the unit normal vector field on~$\Sigma_u$, and express it in terms of the splitting~$\R^n = T_p\Sigma \oplus \operatorname{span}\{\nu\}$ as
    \begin{equation}\label{EQ::nuTildeForm}
        \tilde{\nu}(\Phi_u) = \frac{\nu + w}{\sqrt{1+g_{ij} w^i w^j}}
    \end{equation}
    for some~$w = w^j\partial_j \in T_p\Sigma$. The space~$T_{\Phi_u(p)}\Sigma_u$ is spanned by the vectors
    \begin{equation}\label{EQ::TphifpBasis}
        E_k \coloneq d\Phi_u[\partial_k] = \partial_k + uh_k^l \partial_l + \pd{u}{k} \nu.
    \end{equation}
    By definition,~$\tilde{\nu}$ must be such that, for every~$k = 1,\dots,n-1$,
    \begin{equation}\label{EQ::nuTildeNormalCondition}
        \left\langle E_k, \tilde{\nu}(\Phi_u) \right\rangle = 0.
    \end{equation}
    Substituting~\eqref{EQ::nuTildeForm} and~\eqref{EQ::TphifpBasis} into the orthogonality condition~\eqref{EQ::nuTildeNormalCondition}, we have that
    \begin{equation}\label{EQ::nuTildeCoefficientsComputation}
        \begin{split}
            0&=\left\langle \partial_k + uh_k^l \partial_l + \pd{u}{k} \nu, \frac{\nu + w^m\partial_m}{\sqrt{1+g_{ij} w^i w^j }} \right\rangle\\
            &= \frac{1}{\sqrt{1+g_{ij} w^i w^j}}\bigg(w^m\langle \partial_k, \partial_m\rangle + u h_k^l w^m\langle \partial_l, \partial_m \rangle + \pd{u}{k}\langle\nu,\nu\rangle\bigg)\\
            &= \frac{1}{\sqrt{1+g_{ij} w^i w^j}}\bigg(\big(\delta_k^l + uh_k^l\big)w^m g_{ml} + \pd{u}{k}\bigg),
        \end{split}
    \end{equation}
    where we also used the orthogonality between~$\nu$ and~$\partial_k$.
    
    Since~$u \in \mathcal{F}(\Sigma)$, the term~$\big(\delta_k^l + uh_k^l\big)$ is everywhere non-degenerate\footnote{In particular,~$\rho(\Sigma)$ is bounded by the infimum of the reciprocals of the principal curvatures of~$\Sigma$, which coincide with the eigenvalues of the quadratic form associated with~$h_k^l$. Since, by definition, every~$u \in \mathcal{F}(\Sigma)$ is pointwise bounded by~$\frac{\rho(\Sigma)}{2}$, all the eigenvalues of~$uh_k^l$ have absolute value strictly less than~$1$.}. Taking note that the quantity~$\frac{1}{\sqrt{1+g_{ij} w^i w^j}}$ cannot vanish, the system (as~$k$ varies) of equations~\eqref{EQ::nuTildeCoefficientsComputation} uniquely determines (after reindexing) the coefficients~$w^j$ as
    \begin{equation}\label{EQ::nuTildeCoefficients}
        w^j \coloneq -g^{ij}{\left(\delta_i^k + uh_i^k\right)}^{-1} \pd{u}{k} = -g^{kj}\pd{u}{k} + \O{2}{u}.
    \end{equation}

    From the explicit form of the coefficients~\eqref{EQ::nuTildeCoefficients}, we can compute
    \begin{equation}\label{EQ::WDenominatorExplicit}
        \sqrt{1+g_{ij} w^i w^j} = \sqrt{1 + g_{ij} g^{ki}{\left(\delta_k^l + uh_k^l\right)}^{-1}\pd{u}{l} g^{mj}{\left(\delta_m^q + uh_m^q\right)}^{-1} \pd{u}{q}} = 1 + \O{2}{u}.
    \end{equation}
    Substituting~\eqref{EQ::nuTildeCoefficients} and~\eqref{EQ::WDenominatorExplicit} into~\eqref{EQ::nuTildeForm}, the expansion of~$\tilde{\nu}$ becomes
    \begin{equation}\label{EQ::nuTildeExplicit}
        \tilde{\nu}(\Phi_u) = \frac{\nu - g^{ij}{\left(\delta_i^k + uh_i^k\right)}^{-1} \pd{u}{k}\partial_j}{1 + \O{2}{u}} = \nu - g^{kj}\pd{u}{k}\partial_j + \O{2}{u}.
    \end{equation}

    Differentiating the basis elements~$E_j$ defined by~\eqref{EQ::TphifpBasis} we find
    \begin{equation}\label{EQ::EjDerivative}
        \partial_i E_j = \partial_i \big(\delta_j^k + uh_j^k\big) \partial_k+ \big(\delta_j^k + uh_j^k\big) \partial_i\partial_k + \pdd{u}{i}{j} \nu + \pd{u}{j} h_{i}^l\partial_l.
    \end{equation}
    We apply the Gauss formula to the vectors~$\partial_i$, finding
    \begin{equation}\label{EQ::GaussFormula}
        \partial_i\partial_k = \Gamma_{ik}^l\partial_l - h_{ik} \nu.
    \end{equation}
    Moreover,
    \begin{equation}\label{EQ::EjProductRule}
        \partial_i \big(\delta_j^k + uh_j^k\big) = {\pd{u}{i}h_j^k + u \pd{h_j^k}{i}} = \O{1}{u}.
    \end{equation}
    Using~\eqref{EQ::GaussFormula},~\eqref{EQ::EjProductRule}, and~$\pd{u}{j} h_{i}^l\partial_l = \O{1}{u} \partial_l$ we can further simplify the expression in~\eqref{EQ::EjDerivative}, obtaining
    \begin{equation}\label{EQ::EjDerivativeExplicit}
        \partial_i E_j = \bigg(\Gamma_{ij}^l + \O{1}{u}\bigg)\partial_l + \bigg( -h_{ij} - uh_j^k h_{ik} + \pdd{u}{i}{j}\bigg) \nu.
    \end{equation}

    Let us now compute the second fundamental form using the characterization~$\tilde{h}_{ij}(\Phi_u) = -\langle \partial_i E_j, \tilde{\nu} \rangle$, exploiting the expansions~\eqref{EQ::nuTildeExplicit} and~\eqref{EQ::EjDerivativeExplicit} to deduce that
    \begin{equation}\label{EQ::SecondFundamentalFormComputation}
        \begin{split}
            {(\Phi_u^\ast\tilde{h})}_{ij} &= \left\langle\bigg(\Gamma_{ij}^l + \O{1}{u}\bigg)\partial_l, g^{km} \pd{u}{k}\partial_m + \O{2}{u}\right\rangle - \left\langle \bigg( -h_{ij}- uh_j^k h_{ik} + \pdd{u}{i}{j}\bigg) \nu, \nu\right\rangle\\
                &= {\Gamma_{ij}^l g^{km}g_{lm} \pd{u}{k} + h_{ij} + u h_{j}^k h_{ik} - \pdd{u}{i}{j} + g^{km}g_{lm} \pd{u}{k} \O{1}{u}} + \O{2}{u}\\
                &= {\Gamma_{ij}^k \pd{u}{k} + h_{ij} + u h_{j}^k h_{ik} - \pdd{u}{i}{j}} + \O{2}{u}.
        \end{split}
    \end{equation}
    The proof is concluded by substituting the covariant Hessian definition, i.e.~$\hess{u}{i}{j} = \pdd{u}{i}{j} -\Gamma_{ij}^k \pd{u}{k}$, into~\eqref{EQ::SecondFundamentalFormComputation}.
\end{proof}

The computations that we carried out so far allow us to compute the mean curvature of~$\Sigma_u$ at every point in terms of~$u$ and geometric quantities related to~$\Sigma$. In this regard, it is useful to consider the mean curvature operator~$\mathcal{H} \colon \mathcal{F}(\Sigma) \to C^{0,\alpha}(\Sigma)$, defined by its action on a function~$u \in \mathcal{F}$ as
\begin{equation}\label{2.20BIS}
    \mathcal{H}(u)(p) := H[\Sigma_u](\Phi_u(p))\qquad\text{for every }p \in \Sigma.
\end{equation}
The operator~$\mathcal{H}$ is nonlinear, but it enjoys a useful linear expansion given by the next result.

\begin{proposition}\label{PROP::CurvatureExpansion}
    Let~$\Sigma$ be a~$C^{3,\alpha}$, compact, orientable hypersurface embedded in~$\R^n$ and~$u \in \mathcal{F}(\Sigma)$.

    Then, at every point~$p \in \Sigma$, the following expansion holds
    \begin{equation}\label{EQ::CurvatureExpansion}
        \mathcal{H}(u) = H[\Sigma] - \Delta_{\Sigma} u - |A|^2u + \O{2}{u},
    \end{equation}
    where~$\Delta_{\Sigma}$ is the Laplace-Beltrami operator on~$\Sigma$, and~$A$ is the second fundamental form.
\end{proposition}
\begin{proof}
    Let~$x \coloneq (x^1,\dots,x^{n-1})$ be local coordinates around~$p \in U \subset \Sigma$. We compute the mean curvature as the trace of the second fundamental form, i.e.
    \begin{equation}\label{EQ::CurvatureCharacterization}
        \mathcal{H}(u) = H[\Sigma_u](\Phi_u) = \tilde{g}^{ij}(\Phi_u) \tilde{h}_{ij}(\Phi_u) = {(\Phi_u^\ast \tilde{g})}^{ij}{(\Phi_u^\ast \tilde{h})}_{ij}.
    \end{equation}
    From Lemma~\ref{LEM::gTildeExpansion} we infer that
    \begin{equation}\label{EQ::gTildeInverseExpansion}
        {(\Phi_u^\ast \tilde{g})}^{ij} = g^{ij} - 2uh^{ij} + \O{2}{u},
    \end{equation}
    therefore, also employing Lemma~\ref{LEM::hTildeExpansion}, we substitute~\eqref{EQ::hTildeExpansion} and~\eqref{EQ::gTildeInverseExpansion} into~\eqref{EQ::CurvatureCharacterization} to conclude that
    \begin{equation}\label{EQ::CurvatureComputation}
        \begin{split}
            \mathcal{H}(u) &= \bigg(g^{ij} - 2uh^{ij} + \O{2}{u}\bigg)\bigg(h_{ij} - \hess{u}{i}{j} + u h_i^k h_{kj} + \O{2}{u}\bigg)\\
                &= g^{ij}h_{ij} {- g^{ij}\hess{u}{i}{j} + g^{ij}u h_i^k h_{kj} - 2uh^{ij}h_{ij}} + \O{2}{u}\\
                &= H[\Sigma] {- \Delta_{\Sigma}u + u  h_i^k h_{k}^i - 2uh^{ij}h_{ij}} + \O{2}{u}.
        \end{split}
    \end{equation}
    Since~$h_i^k h_{k}^i = h^{ij}h_{ij} = |A|^2$, the result follows from~\eqref{EQ::CurvatureComputation}.
\end{proof}

Proposition~\ref{PROP::CurvatureExpansion} can be seen as a quantitative version of the second variation formula for minimal surfaces. Indeed, as briefly discussed in the introduction, the Jacobi operator~$\mathcal{L}_\Sigma$ defined by~\eqref{EQ::JacobiOperatorDefinition} is strictly related with the stability of a minimal surface. From~\eqref{EQ::CurvatureExpansion} it follows that
\begin{equation}
    \mathcal{H}(u) = H[\Sigma] +\mathcal{L}_\Sigma u + \O{2}{u}.
\end{equation}
Importantly, the nonlinear remainder that appears in the expansion is a term of the type~$\O{2}{u}$, and therefore it is uniformly bounded and uniformly Lipschitz continuous on~$\mathcal{F}(\Sigma)$. This fact will be crucial in the forthcoming proof of Theorem~\ref{THM::Foliation}.

\subsection{Fixed point solutions of the mean curvature equation}
In the following, we exploit the quantitative information offered by Proposition~\ref{PROP::CurvatureExpansion} to treat the nonlinear curvature operator as a small perturbation of the Jacobi operator. In this regard, it is useful to keep in mind that the linearization offered by the expansion~\eqref{EQ::CurvatureExpansion} is slightly different from the one obtained by considering local graphicality over a tangent plane. The advantage, here, is that the deformation is treated globally, allowing for a functional-analytic approach.

We first establish solvability of the linearized problem. The following lemma is a consequence of classical results, but, for the facility of the reader, we state precisely some notable properties that we will use in the rest of this section.

\begin{lemma}\label{LEM::JacobiOperatorProperties}
    Let~$\alpha \in (0,1)$ and let~$\Sigma \subset B_2$ be a compact, orientable,~$C^{3,\alpha}$ hypersurface embedded in~$\R^n$, with boundary~$\partial \Sigma = \Sigma \cap \partial B_2 \neq \varnothing$. Suppose that~$\Sigma$ is a strictly stable minimal surface in~$B_2$. Assume that the~$C^{3,\alpha}$ norm of~$\Sigma$ is bounded by~$M > 0$.

    Then, for every~$f \in C^{0,\alpha}(\Sigma)$ there exists a unique~$u \in C^{2,\alpha}_0(\Sigma)$ such that
    \begin{equation}\label{EQ::JacobiEquationInhomogeneous}
        \mathcal{L}_\Sigma u = f.
    \end{equation}

    Moreover, there exists a constant~$C > 0$, which depends only on~$n$,~$\alpha$,~$\lambda_1(\Sigma)$,~$M$, and~$\rho_{\loc}(\Sigma)$, such that
    \begin{equation}\label{EQ::JacobiSchauder}
        \|u\|_{C^{2,\alpha}(\Sigma)} \leq C\|f\|_{C^{0,\alpha}(\Sigma)}.
    \end{equation}

    Furthermore, if~$f \geq 0$ on~$\Sigma$, then~$u \geq 0$ on~$\Sigma$.
\end{lemma}
\begin{proof}
    The operator~$\mathcal{L}_\Sigma$ is linear and local. Therefore, if~$\Sigma$ is not connected we can solve the equation~\eqref{EQ::JacobiEquationInhomogeneous} by finding~$u^j \in \mathcal{F}(\Sigma^j)$ such that
    \begin{equation}
        \mathcal{L}_{\Sigma^j} u^j = \restr{f}{\Sigma^j}
    \end{equation}
    and defining~$u \coloneq \sum_j u^j$. In particular, since
    \begin{equation}
        \lambda_1(\Sigma) = \inf_j \lambda_1(\Sigma^j)
    \end{equation}
    and using~\eqref{EQ::rhoLocDefinition}, the constant~$C$ in~\eqref{EQ::JacobiSchauder} can be taken as the maximum among the constants on each~$\Sigma^j$. Hence, from now on we will assume that~$\Sigma$ is connected.

    Consider the bilinear form~$\mathcal{E} \colon H^1_0(\Sigma) \times H^1_0(\Sigma) \to \R$ defined by
    \begin{equation}
        \mathcal{E}(u,v) \coloneq \int_{\Sigma} \left(\langle\nabla_\Sigma u, \nabla_\Sigma v\rangle - |A|^2 uv\right)\,d\mathcal{H}^{n-1},
    \end{equation}
    where~$\nabla_\Sigma$ is the tangential gradient. By applying the divergence theorem, we see that every distributional solution~$u$ of~\eqref{EQ::JacobiEquationInhomogeneous} is such that
    \begin{equation}
        \mathcal{E}(u,v) = F(v)\coloneq \int_\Sigma fv\,d\mathcal{H}^{n-1}
    \end{equation}
    for every~$v \in H^1_0(\Sigma)$. Applying a density argument to~\eqref{EQ::StabilitySpectralDefinition} we see that
    \begin{equation}\label{EQ::CoercivityL2}
        \mathcal{E}(u,u) \geq \lambda_1(\Sigma) \|u\|^2_{L^2(\Sigma)},
    \end{equation}
    thus, using that~$|A|^2 \leq (n-1)M^2$ everywhere, we have that
    \begin{equation}
        \begin{split}
            \|\nabla_\Sigma u\|^2_{L^2(\Sigma)} &= \mathcal{E}(u,u) + \int_\Sigma |A|^2 u^2\,d\mathcal{H}^{n-1} \\
                &\leq \mathcal{E}(u,u) + (n-1)M^2 \|u\|^2_{L^2(\Sigma)}\\
                &\leq \left(1+\frac{(n-1)M^2}{\lambda_1(\Sigma)}\right)\mathcal{E}(u,u),
        \end{split}
    \end{equation}
    proving the coercivity of~$\mathcal{E}$. Similarly,
    \begin{equation}
        \begin{split}
            \mathcal{E}(u,v) &\leq \int_\Sigma \langle\nabla_\Sigma u, \nabla_\Sigma v\rangle\,d\mathcal{H}^{n-1} + \int_{\Sigma}|A|^2 uv\,d\mathcal{H}^{n-1} \\
            &\leq \|\nabla_\Sigma u\|_{L^2(\Sigma)}\|\nabla_\Sigma v\|_{L^2(\Sigma)} + (n-1)M^2 \|u\|_{L^2(\Sigma)}\|v\|_{L^2(\Sigma)}
        \end{split}
    \end{equation}
    and
    \begin{equation}
        F(v) \leq \|f\|_{L^{\infty}(\Sigma)} \|v\|_{L^2(\Sigma)}.
    \end{equation}
    All in all, the Lax-Milgram theorem (see e.g.~\cite{MR2759829}) yields a weak solution~$u \in H^1_0(\Sigma)$ for the equation~\eqref{EQ::JacobiEquationInhomogeneous}. Note also that
    \begin{equation}\label{EQ::JacobiLaxMilgramStability}
        \|u\|_{H^1(\Sigma)} \leq C(n,M,\lambda_1(\Sigma))\|f\|_{L^{\infty}(\Sigma)}.
    \end{equation}

    To conclude the proof of~\eqref{EQ::JacobiSchauder}, we apply a standard bootstrap argument by covering~$\Sigma$ with local charts. Since each flattening map (together with its inverse) is a diffeomorphism with~$C^{3,\alpha}$ norm bounded only in terms of~$M$, the regularity constant in~\eqref{EQ::JacobiSchauder} will depend only on~$\alpha$,~$M$, and the constant~$C$ that appears in~\eqref{EQ::JacobiLaxMilgramStability}.

    Finally, we prove the weak maximum principle for the operator~$\mathcal{L}_\Sigma$. Suppose that~$f \geq 0$ and let~$w \coloneq u_- = \max{\{-u, 0\}}$. Leveraging~\eqref{EQ::CoercivityL2} and the fact that~$w$ and~$-u$ coincide on the support of~$w$, we have that
    \begin{equation}
        0 \leq F(w) \leq \mathcal{E}(u,w) = -\mathcal{E}(w,w) \leq -\lambda_1(\Sigma)\|w\|_{L^2(\Sigma)}^2 \leq 0,
    \end{equation}
    where we also used the non-negativity of~$f$ and~$w$, and the fact that~$u$ is a solution. Therefore, since~$\lambda_1(\Sigma) > 0$, we find that~$w = 0$ almost everywhere on~$\Sigma$, and thus~$u \geq 0$ everywhere since it is continuous.
\end{proof}

Let us now consider the nonlinear remainder~$\mathcal{Q} \colon \mathcal{F}(\Sigma) \to C^{0,\alpha}(\Sigma)$, defined by
\begin{equation}\label{EQ::QuDefinition}
    \mathcal{Q} \colon u \longmapsto \mathcal{H}(u) - \mathcal{L}_\Sigma u.
\end{equation}
Thanks to Proposition~\ref{PROP::CurvatureExpansion}, if~$\Sigma$ has zero mean curvature, then
\begin{equation}\label{EQ::QisO2u}
    \mathcal{Q}(u) = \O{2}{u}.
\end{equation}
Let~$t \in \R$ and consider the operator family~$\mathcal{A}_t \colon \mathcal{F}(\Sigma) \to C^{2,\alpha}_0(\Sigma)$, defined by
\begin{equation}\label{EQ::AtOperatorDefinition}
    \mathcal{A}_t \colon u \longmapsto \mathcal{L}_\Sigma^{-1}(t - \mathcal{Q}(u)).
\end{equation}
Thanks to Lemma~\ref{LEM::JacobiOperatorProperties} and~\eqref{EQ::QisO2u}, the operator~$\mathcal{A}_t$ is well defined whenever~$\Sigma$ is a strictly stable minimal surface. The role of the operator~$\mathcal{A}_t$ is clarified by the following result.

\begin{lemma}\label{LEM::FixedPointIsCMC}
    Let~$\alpha \in (0,1)$ and let~$\Sigma \subset B_2$ be a compact, orientable,~$C^{3,\alpha}$ hypersurface embedded in~$\R^n$, with boundary~$\partial \Sigma = \Sigma \cap \partial B_2 \neq \varnothing$. Suppose that~$\Sigma$ is a strictly stable minimal surface in~$B_2$. Assume that the~$C^{3,\alpha}$ norm of~$\Sigma$ is bounded by~$M > 0$. Let~$t \in \R$ and~$\mathcal{A}_t$ be defined by~\eqref{EQ::AtOperatorDefinition}.

    Suppose that~$u \in \mathcal{F}(\Sigma)$ is a fixed point of the operator~$\mathcal{A}_t$, i.e.,
    \begin{equation}
        u = \mathcal{A}_t(u).
    \end{equation}

    Then, the hypersurface~$\Sigma_u$, as defined in~\eqref{EQ::SigmaUDefinition}, is a constant mean curvature hypersurface, with
    \begin{equation}\label{EQ::CMCforFixedPoint}
        H[\Sigma_u] \equiv t.
    \end{equation}
\end{lemma}
\begin{proof}
    Using Proposition~\ref{PROP::CurvatureExpansion} and~\eqref{EQ::QisO2u}, we can compute
    \begin{equation}
        H[\Sigma_u] = \mathcal{H}(u) = {\mathcal{L}_{\Sigma} u} + \mathcal{Q}(u),
    \end{equation}
    hence, using that~$u = \mathcal{A}_t(u)$, we find
    \begin{equation}
        H[\Sigma_u] = {\mathcal{L}_{\Sigma} \mathcal{A}_t(u)} + \mathcal{Q}(u) = {\mathcal{L}_{\Sigma} \mathcal{L}_{\Sigma}^{-1} (t - \mathcal{Q}(u))} + \mathcal{Q}(u) = t.\qedhere{}
    \end{equation}
\end{proof}
It is necessary to add some assumptions on the value~$t \in \R$ in order for the operator~$\mathcal{A}_t$ to admit a fixed point, see also Appendix~\ref{APP::PlaneExample}. Our goal is to show that, if~$|t|$ is quantitatively small enough, a fixed point always exists. This is possible by estimating the boundedness and Lipschitz constant of~$\mathcal{A}_t$ in an appropriate closed subset of~$C^{2,\alpha}_0(\Sigma)$. The next result shows that the restriction of~$\mathcal{A}_t$ to such closed set is a contraction, and thus has a fixed point.
\begin{lemma}\label{LEM::FixedPoint}
    Let~$\alpha \in (0,1)$ and let~$\Sigma \subset B_2$ be a compact, orientable,~$C^{3,\alpha}$ hypersurface embedded in~$\R^n$, with boundary~$\partial \Sigma = \Sigma \cap \partial B_2 \neq \varnothing$. Suppose that~$\Sigma$ is a strictly stable minimal surface in~$B_2$. Assume that the~$C^{3,\alpha}$ norm of~$\Sigma$ is bounded by~$M > 0$. Let~$t \in \R$ and~$\mathcal{A}_t$ be defined by~\eqref{EQ::AtOperatorDefinition}.

    Then, there exists a constant~$c > 0$, which depends only on~$n$,~$\alpha$,~$\lambda_1(\Sigma)$,~$M$, and~$\rho_{\loc}(\Sigma)$, such that, if
    \begin{equation}\label{EQ::tSmallness}
        |t| \leq c \min\{1,\sep(\Sigma)\},
    \end{equation}
    the operator~$\mathcal{A}_t$ admits a fixed point~$u_t \in \mathcal{F}(\Sigma)$.

    Moreover, there exists a constant~$C_\ast > 0$, which depends only on~$n$,~$\alpha$,~$\lambda_1(\Sigma)$,~$M$, and~$\rho_{\loc}(\Sigma)$, such that
    \begin{equation}\label{EQ::uniformPerturbationBound}
        \|u_t\|_{C^{2,\alpha}(\Sigma)} \leq C_\ast.
    \end{equation}
\end{lemma}
\begin{proof}
    Let~$C_0$ be the constant in~\eqref{EQ::JacobiSchauder} given by Lemma~\ref{LEM::JacobiOperatorProperties}.
    Recalling~\eqref{EQ::QisO2u}, let~$C_1$ be the constant for which~$\mathcal{Q}$ satisfies~\eqref{EQ::OkNotationDefinition}. For any~$u \in \mathcal{F}(\Sigma)$, we estimate
    \begin{equation}\label{EQ::AtUSizeEstimate}
        \|\mathcal{A}_t(u)\|_{C^{2,\alpha}(\Sigma)} = \|\mathcal{L}_\Sigma^{-1}({t - \mathcal{Q}(u)})\|_{C^{2,\alpha}(\Sigma)} \leq C_0\left(\|\mathcal{Q}(u)\|_{C^{0,\alpha}(\Sigma)} + \|t\|_{C^{0,\alpha}(\Sigma)}\right).
    \end{equation}
    Clearly,~$\|t\|_{C^{0,\alpha}(\Sigma)} = |t|$. Therefore,~\eqref{EQ::AtUSizeEstimate} becomes
    \begin{equation}\label{EQ::AtUSizeEstimateRefined}
        \|\mathcal{A}_t(u)\|_{C^{2,\alpha}(\Sigma)} \leq C_0 C_1 \|u\|_{C^{2,\alpha}(\Sigma)}^2 + C_0 |t|.
    \end{equation}
 Similarly, for~$u,v \in \mathcal{F}(\Sigma)$, the map~$\mathcal{A}_t$ satisfies an almost-Lipschitz condition
    \begin{equation}\label{EQ::almostLipschitz}
        \begin{split}
            \|\mathcal{A}_t(u) - \mathcal{A}_t(v)\|_{C^{2,\alpha}(\Sigma)} &= \|\mathcal{L}_\Sigma^{-1}({\mathcal{Q}(v) - \mathcal{Q}(u)})\|_{C^{2,\alpha}(\Sigma)} \\
                &\leq C_0 C_1 (\|u\|_{C^{2,\alpha}(\Sigma)} + \|v\|_{C^{2,\alpha}(\Sigma)}) \|u - v\|_{C^{2,\alpha}(\Sigma)}.
        \end{split}
    \end{equation}

    Let
    \begin{equation}\label{EQ::FixedPointRadiusDefinition}
        r \coloneq \min\left\{\frac{1}{8C_0C_1}, \frac{1}{3}\rho(\Sigma)\right\}
    \end{equation}
    and
    \begin{equation}
        \mathcal{B} \coloneq \{u \in C^{2,\alpha}_0(\Sigma) \colon \|u\|_{C^{2,\alpha}(\Sigma)} \leq r\}.
    \end{equation}
    The set~$\mathcal{B}$ is a closed subset of~$C^{2,\alpha}_0(\Sigma)$ and, thanks to~\eqref{EQ::FixedPointRadiusDefinition} we have that~$r \leq \frac{1}{3}\rho(\Sigma)$ and therefore~$\mathcal{B} \subset \mathcal{F}(\Sigma)$. Since~$C_0C_1 r \leq \frac{1}{8}$, for any~$u \in \mathcal{B}$ we can take advantage of~\eqref{EQ::AtUSizeEstimateRefined} to conclude that
    \begin{equation}
        \|\mathcal{A}_t(u)\|_{C^{2,\alpha}(\Sigma)} \leq C_0C_1 r^2 + C_0|t| \leq \frac{r}{8} + C_0 |t|.
    \end{equation}
    From this, if we additionally suppose that
    \begin{equation}\label{EQ::tSmallAbsValue}
        |t| \leq \frac{7}{8C_0}r,
    \end{equation}
    we deduce that
    \begin{equation}
        \|\mathcal{A}_t(u)\|_{C^{2,\alpha}(\Sigma)} \leq r,
    \end{equation}
    hence,
    \begin{equation}\label{EQ::AtSelfMapOnB}
        \mathcal{A}_t(\mathcal{B}) \subset \mathcal{B}.
    \end{equation}
    In addition, specializing~\eqref{EQ::almostLipschitz} to~$u,v \in \mathcal{B}$ yields
    \begin{equation}\label{EQ::AtLipschitzOnB}
        \|\mathcal{A}_t(u) - \mathcal{A}_t(v)\|_{C^{2,\alpha}(\Sigma)} \leq 2C_0C_1 r \|u - v\|_{C^{2,\alpha}(\Sigma)} \leq \frac{1}{4}\|u - v\|_{C^{2,\alpha}(\Sigma)}.
    \end{equation}

    Gathering~\eqref{EQ::AtSelfMapOnB} and~\eqref{EQ::AtLipschitzOnB}, the map~$\mathcal{A}_t$ is a contraction on~$\mathcal{B}$ if~$t$ satisfies~\eqref{EQ::tSmallAbsValue}. On this account, we let
    \begin{equation}\label{EQ::smallcDefintion}
        c \coloneq \min\left\{\frac{7}{64C_0^2C_1},\frac{\min\{1,\rho_{\loc}(\Sigma)\}}{4C_0}\right\}
    \end{equation}
    and we recall that~$\rho(\Sigma) = \min\{\sep(\Sigma),\rho_{\loc}(\Sigma)\}$. Therefore, also using~\eqref{EQ::FixedPointRadiusDefinition}, the condition in~\eqref{EQ::tSmallAbsValue} is satisfied as a consequence of~\eqref{EQ::tSmallness} with~$c$ as in~\eqref{EQ::smallcDefintion}. Indeed,
    \begin{equation}
        |t| \leq c \leq \frac{7}{64C_0^2C_1} = \frac{7}{8C_0} \frac{1}{8C_0C_1}
    \end{equation}
    and
    \begin{equation}
        |t| \leq c\min\{1,\sep(\Sigma)\} \leq \min\{1,\sep(\Sigma)\}\frac{\min\{1,\rho_{\loc}(\Sigma)\}}{4C_0} \leq \frac{\rho(\Sigma)}{4C_0} \leq \frac{7}{8C_0}\frac{1}{3}\rho(\Sigma).
    \end{equation}
    The contraction property of~$\mathcal{A}_t$ allows us to apply the Banach-Caccioppoli Theorem to find the desired fixed point~$u_t$.

    Since~$u_t \in \mathcal{B}$, we can choose
    \[C_\ast \coloneq \frac{1}{8C_0C_1} \geq r \geq \|u\|_{C^{2,\alpha(\Sigma)}}\]
    in~\eqref{EQ::uniformPerturbationBound}, concluding the proof.
\end{proof}

We want to emphasize that the choice of~$r$ in~\eqref{EQ::FixedPointRadiusDefinition} achieves two different goals. On one hand, it ensures that the map~$\mathcal{A}_t$ is a contraction on the set~$\mathcal{B}$. On the other, we need all the perturbed components to avoid touching each other. This translates directly to the definition of~$c$ in~\eqref{EQ::smallcDefintion}, shedding light on the condition~\eqref{EQ::eps0OptimalCondition}.

Roughly speaking, the geometry of each connected component of~$\Sigma$ naturally limits the maximal mean curvature that can be achieved by perturbations supported in~$B_1$. However, when multiple components are considered, the only obstacle to the realization of the optimal perturbation of each component is exactly the limitation we impose by asking the perturbed surface to avoid self-intersections.

\subsection{Conclusion of the Proof of Theorem~\ref{THM::Foliation}}
    Let~$c$ be as in Lemma~\ref{LEM::FixedPoint} and
    \begin{equation}
        \eps_0 \coloneq c \min\{1, \sep(\Sigma)\}.
    \end{equation}
    For~$t \in (-\eps_0, \eps_0)$, let~$u_t$ be the fixed point given by Lemma~\ref{LEM::FixedPoint}.
    
    We claim that the family~${\{\Sigma_t \coloneq \Sigma_{u_t}\}}_{t \in (-\eps_0,\eps_0)}$ constitutes the desired foliation. First, we notice that~\eqref{EQ::CMCforFixedPoint} in Lemma~\ref{LEM::FixedPointIsCMC} immediately gives~\eqref{EQ::FoliationCMCCondition}. Therefore, it only remains to prove that~${\{\Sigma_t\}}_{t}$ is ordered and its union covers a neighborhood of~$\Sigma \smallsetminus \partial \Sigma$.

    Let~$t,\tau \in (-\eps_0,\eps_0)$ with~$t < \tau$. Using~\eqref{EQ::FoliationCMCCondition},
    and recalling the notation in~\eqref{2.20BIS}, we have that
    \begin{equation}
        \mathcal{H}(u_\tau) - \mathcal{H}(u_t) = \tau - t.
    \end{equation}
    From this and~\eqref{EQ::QuDefinition}, it follows that
    \begin{equation}\label{EQ::tTauDifferencePDE}
        \mathcal{L}_\Sigma({u_\tau - u_t}) = \tau - t + \mathcal{Q}(u_t) - \mathcal{Q}(u_\tau).
    \end{equation}
    As in the proof of Lemma~\ref{LEM::FixedPoint}, we denote by~$C_0$ the constant in~\eqref{EQ::JacobiSchauder} given by Lemma~\ref{LEM::JacobiOperatorProperties}, and, recalling~\eqref{EQ::QisO2u}, by~$C_1$ the constant for which~$\mathcal{Q}$ satisfies~\eqref{EQ::OkNotationDefinition}. Using~\eqref{EQ::OkNotationDefinition} we obtain
    \begin{equation}\label{EQ::ututauRemainderEstimate}
        \|\mathcal{Q}(u_\tau) - \mathcal{Q}(u_t)\|_{C^{0,\alpha}(\Sigma)} \leq 2C_1 r\|u_\tau - u_t\|_{C^{2,\alpha}(\Sigma)},
    \end{equation}
    where~$r$ is given in~\eqref{EQ::FixedPointRadiusDefinition}.
    
    Thus, applying Lemma~\ref{LEM::JacobiOperatorProperties} with~$u \coloneq {u_\tau - u_t}$, we find that
    \begin{equation}\label{EQ::ututauBeforeAbsorb}
        \|u_\tau - u_t\|_{C^{2,\alpha}(\Sigma)} \leq C_0 (\tau - t) + 2C_0C_1 r\|u_\tau - u_t\|_{C^{2,\alpha}(\Sigma)}.
    \end{equation}
    Hence, using that~$2C_0C_1 r \leq \frac{1}{4}$, we absorb the term in the right hand side of~\eqref{EQ::ututauBeforeAbsorb}, which in turn becomes
    \begin{equation}\label{EQ::utContinuous}
        \|u_\tau - u_t\|_{C^{2,\alpha}(\Sigma)} \leq \frac{4}{3}C_0(\tau - t). 
    \end{equation}
    Plugging this into~\eqref{EQ::ututauRemainderEstimate}, we conclude that
    \begin{equation}
        \|\mathcal{Q}(u_t) - \mathcal{Q}(u_\tau)\|_{L^{\infty}(\Sigma)} \leq \|\mathcal{Q}(u_t) - \mathcal{Q}(u_\tau)\|_{C^{0,\alpha}(\Sigma)} \leq \frac{8}{3}C_0C_1 r (\tau - t) \leq \frac{1}{3} (\tau - t).
    \end{equation}
    
    In particular, the right hand side in~\eqref{EQ::tTauDifferencePDE} is strictly positive everywhere on~$\Sigma$. Therefore, Lemma~\ref{LEM::JacobiOperatorProperties} entails that~${u_\tau - u_t}$ is non-negative on~$\Sigma$, and therefore the surfaces~$\Sigma_{u_\tau}$ and~$\Sigma_{u_t}$ are weakly ordered. Since they share the same boundary and have different curvature, the only possibility is that their relative interiors are disjoint, yielding strong ordering of~${\{\Sigma_t\}}_{t \in (-\eps_0,\eps_0)}$.

\begin{figure}[ht]
        \centering
        \includegraphics[width=0.8\textwidth]{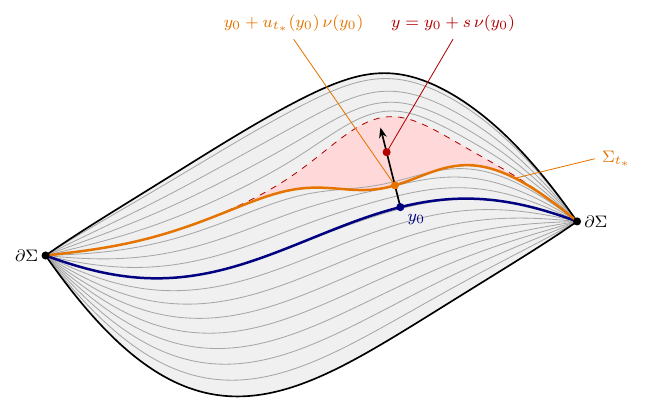}
        \caption{Sketch illustrating the `hole' not reached by the foliation, inside which~$y$ lies. Here, we depict the impossible situation in which~$u_{t_\ast}(y_0) < s$.}\label{FIG::Hole}
    \end{figure}

    Consider the set~$U \subset \R^n$ such that~$\partial U = \Sigma_{\eps_0} \cup \Sigma_{-\eps_0}$ and~$U$ is open and bounded. Suppose now by contradiction that~${\{\Sigma_t\}}_{t \in (-\eps_0,\eps_0)}$ fails to cover a neighborhood of~$\Sigma \smallsetminus \partial \Sigma$. In other words, suppose that there exists a point
    \begin{equation}
        y \in U \smallsetminus \bigcup_{t \in (-\eps_0, \eps_0)} \Sigma_t.
    \end{equation}
    Notice that the set~$U$ is contained in a~$r$-neighborhood of~$\Sigma$. As a consequence,~$y$ admits a unique projection~$y_0 \in \Sigma$ such that
    \begin{equation}\label{EQ::yProjection}
        y = y_0 + s \nu(y_0),
    \end{equation}
    where~$s \in (-r,r)$ is the signed distance of~$y$ from~$\Sigma$. Let
    \begin{equation}\label{EQ::tastDefinition}
        t_\ast \coloneq \sup\{t \in (-\eps_0, \eps_0) \colon u_t(y_0) \leq s\}.
    \end{equation}
    See also Figure~\ref{FIG::Hole} for a representation of the construction of~$t_\ast$.

    Recalling~\eqref{EQ::utContinuous}, the map~$(-\eps_0,\eps_0) \ni t \mapsto u_t \in C^{2,\alpha}(\Sigma)$ is continuous. Accordingly,
    \begin{equation}
        u_{t_\ast}(y_0) = s.
    \end{equation} 
    This and~\eqref{EQ::yProjection} imply that~$y = y_0 + u_{t_\ast}(y_0) \nu(y_0)$, and therefore~$y \in \Sigma_{t_\ast}$. Although~\eqref{EQ::tastDefinition} does not apriori exclude the possibility that~$t_\ast =\eps_0$, we notice that~$y_0 + u_{\eps_0}(y_0) \nu(y_0) \in \Sigma_{\eps_0} \subset \partial U$, which contradicts the fact that~$y$ belongs to the open set~$U$.
\qed{}

\section{Existence of stable s-minimal surfaces and Proof of Theorem~\ref{THM::Main}}\label{SEC::sMinimalSurface}
In this section we leverage the foliation provided by Theorem~\ref{THM::Foliation} to prove Theorem~\ref{THM::Main} by combining fine estimates for the nonlocal mean curvature and a barrier argument relying on the maximum principle for the nonlocal mean curvature of regular hypersurfaces.

\subsection{Nonlocal barriers and perimeter minimization}
This section contains the construction of a stable minimal surface trapped between two strict classical barriers for a suitable fractional parameter~$s \in (0,1)$. As we will show, one can pick a value of~$s$ that is large enough for classical barriers to act as nonlocal barriers. We stress that throughout this section we never assume connectedness of sets, hence the construction works in the same way for both the cases of~$s$-minimal surfaces of single and multiple connected components in a ball.

We start off with a technical lemma, finding a quantitative uniform graphical radius for a~$C^2$ hypersurface.

\begin{figure}[ht]
\centering
\includegraphics[width=0.5\textwidth]{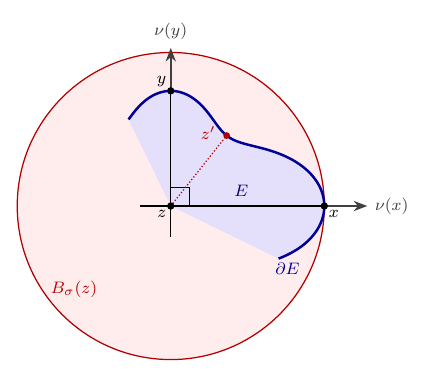}
\caption{Sketch of the contradiction found in the proof of Lemma~\ref{LEM::UniformGraphRadius}. Note that, here,~$r$ is chosen much larger than the reach of~$\partial E$, contrary to the assumptions of Lemma~\ref{LEM::UniformGraphRadius}.}\label{FIG::GraphRadius}
\end{figure}

\begin{lemma}\label{LEM::UniformGraphRadius}
    Let~$E \subset \R^n$ be such that~$\partial E \cap B_2$ is a~$C^2$ hypersurface. Let~$r \coloneq \min\{1,\rho(\partial E \cap B_2)\}$.

    Then, for every~$x \in \partial E \cap B_1$, we have that~$\partial E \cap B_{r}(x)$ is a~$C^2$ graph in some direction.
\end{lemma}
\begin{proof}
    Let~$x \in \partial E \cap B_1$. Clearly,~$B_{r}(x) \subset B_{2}$. Therefore, since
    \begin{equation}\label{EQ::rSmallerThanReach}
        r \leq \rho(\partial E \cap B_{2}),
    \end{equation}
    it follows that~$\partial E \cap B_{r}(x)$ is a connected hypersurface, with all its principal curvatures bounded by~$1/r$. As a consequence, 
    \begin{equation}\label{3.1BIS}\begin{split}&
    {\mbox{at every~$y \in \partial E \cap B_r(x)$, the set~$\partial E$ is touched}}\\&{\mbox{on both sides of its orientation by a ball of radius~$r$.}}\end{split}\end{equation}

    Denote by~$\nu$ the outward unit normal of~$\partial E$. Since~$\partial E \cap B_{r}(x)$ is a connected~$C^2$ surface, to obtain the result we just need to prove that
    \begin{equation}\label{EQ::NormalNotOrthogonal}
        \langle \nu(x), \nu(y) \rangle \neq 0 
    \end{equation}
    for every~$y \in \partial E \cap B_{r}(x)$. In fact, this would entail that~$\partial E \cap B_{r}(x)$ is a graph over the tangent hyperplane at~$x$.
    
    Suppose by contradiction that~$y \in \partial E \cap B_{r}(x)$ is such that~\eqref{EQ::NormalNotOrthogonal} is false (see Figure~\ref{FIG::GraphRadius}). Then, the two lines~$\{x + s\nu(x) \colon s \in \R\}$ and~$\{y + t \nu(y) \colon t \in \R\}$ must meet orthogonally at a point~$z = x+s\nu(x) = y+t\nu(y)$ for some~$s,t \in \R$. By the Pythagorean Theorem,
    \begin{equation}\label{EQ::stPythagoreanTheorem}
        \sigma \coloneq \max\{|x-z|,|y-z|\} = \max\{|s|,|t|\} \leq \sqrt{s^2+t^2} = |x-y| < r.
    \end{equation}
    From this, recalling~\eqref{EQ::rSmallerThanReach}, we have that~$z$ has a unique projection~$z^\prime \in \partial E \cap B_{2}$, which attains the optimal distance, i.e.,~$|z-z^\prime| = \operatorname{dist}(z,\partial E \cap B_{2})$. The uniqueness of~$z^\prime$ entails that~$\sigma > |z-z^\prime|$. Using~\eqref{3.1BIS} and~\eqref{EQ::stPythagoreanTheorem}, we infer that the ball~$B_\sigma(z)$, which is tangent to~$\partial E$ either at the point~$x$ or at the point~$y$, or both, must be entirely contained in either~$E$ or~$E^c$. However, this contradicts the fact that~$z^\prime \in \partial E$ is an interior point of such ball. This contradiction shows~\eqref{EQ::NormalNotOrthogonal} and concludes the proof.
\end{proof}

The next partial result leverages Lemma~\ref{LEM::UniformGraphRadius} to pass from a lower bound on the mean curvature to a lower bound on the~$s$-mean curvature for large~$s$. The precise statement goes as follows.

\begin{lemma}\label{LEM::sMeanCurvatureLowerBound}
    Let~$E \subset \R^n$ be such that~$\partial E \cap B_{2}$ is a~$C^{2,\alpha}$ hypersurface for some~$\alpha \in (0,1)$. Suppose that the~$C^{2,\alpha}$ norm of~$\partial E \cap B_{2}$ is bounded by~$M > 0$. Let~$r \coloneq \min\{1,\rho(\partial E \cap B_{2})\}$. Additionally, assume that there exists~$\lambda > 0$ such that, for every~$x \in \partial E \cap B_{1}$,
    \begin{equation}\label{EQ::MeanCurvatureLowerBound}
        H[\partial E](x) \geq \lambda.
    \end{equation}

    Then, there exists~$s_0 \in (0,1)$ such that, for every~$s \in [s_0,1)$ and~$x \in \partial E \cap B_{1}$, there holds %chktex 9
    \begin{equation}\label{EQ::sMeanCurvatureLowerBound}
        H_s[\partial E](x) \geq 1.
    \end{equation}

    Moreover, there exists a constant~$c > 0$, which depends only on~$n$,~$\alpha$, and~$M$, such that the smallest
    value of~$s_0$ for which the above statement holds true satisfies
    \begin{equation}\label{EQ::optimalS0Estimate}
        s_0 \leq \max\left\{\frac{1}{10},1-c r \lambda\right\}.
    \end{equation}
\end{lemma}
\begin{proof}
    Let~$x \in \partial E \cap B_1$. Thanks to Lemma~\ref{LEM::UniformGraphRadius}, the set~$\partial E \cap B_r(x)$ is a~$C^{2,\alpha}$ graph. Consider the blowup map~$T \colon \R^n \to \R^n$ defined by
    \begin{equation}\label{EQ::TBlowupDefinition}
        T \colon y \longmapsto \frac{y-x}{r}.
    \end{equation}
    Let~$E_\ast \coloneq T(E)$. Since~$T^{-1}(B_1) = B_r(x) \subset B_{2}$ and~$r \le 1$, the set~$\partial E_\ast \cap B_1$ is a~$C^{2,\alpha}$ graph, and its norm is still bounded by~$M$. We apply~\cite{MR3107529}*{Lemma~9} for any~$s \in [1/10,1)$, with~$\alpha \coloneq \alpha$,~$F \coloneq E_\ast$,~$x_0 = 0$, and~$M$, and we find~$c_0 = c_0(n) > 0$ and~$C = C(n,\alpha) > 0$ such that %chktex 9
    \begin{equation}
        \bigg|H[\partial E_\ast](0) - c_0(1-s)H_s[\partial E_\ast](0)\bigg|\leq CM(1-s),
    \end{equation}
    thus,
    \begin{equation}\label{EQ::sMeancurvatureAsymptotic}
        H_s[\partial E_\ast](0) \geq \frac{1}{c_0}\left(\frac{1}{1-s}H[\partial E_\ast](0) - CM\right).
    \end{equation}

    We use~\eqref{EQ::MeanCurvatureLowerBound}, together with translation invariance and scale homogeneity of the mean curvature to conclude that
    \begin{equation}
        H[\partial E_\ast](0) \geq r\lambda.
    \end{equation}
    Plugging this into~\eqref{EQ::sMeancurvatureAsymptotic}, we obtain
    \begin{equation}\label{EQ::sCurvatureLowerBoundWithS}
        H_s[\partial E_\ast](0) \geq \frac{1}{c_0}\left(\frac{r\lambda}{1-s} - CM\right).
    \end{equation}

    Let
    \begin{equation}\label{EQ::optimalS0Value}
        s_0 \coloneq \max\left\{\frac{1}{10}, 1-\frac{r\lambda}{c_0+CM}\right\}.
    \end{equation}
    Then, for every~$s \in [s_0,1)$, we have that %chktex 9
    \begin{equation}
        1-s \leq \frac{r\lambda}{c_0+CM},
    \end{equation}
    hence, we can rewrite~\eqref{EQ::sCurvatureLowerBoundWithS} as~$H_s[\partial E_\ast](0) \geq 1$. From this, and recalling that~$E_\ast = T(E)$ with~$T$ defined by~\eqref{EQ::TBlowupDefinition}, we can exploit the~$s$-homogeneity of the~$s$-mean curvature to find
    \begin{equation}
        H_s[\partial E](x) \geq r^{-s},
    \end{equation}
    which gives~\eqref{EQ::sMeanCurvatureLowerBound} because~$r \leq 1$.

    Clearly,~\eqref{EQ::optimalS0Estimate} immediately follows from~\eqref{EQ::optimalS0Value} upon choosing
    \[c \coloneq \frac{1}{c_0 + CM}.\qedhere{}\]
\end{proof}

\begin{remark}\label{REM::sMeanCurvatureLowerBoundSign}
    We remark that the mean curvature is defined up to a sign change. Indeed, by reversing the orientation of a hypersurface, which is equivalent to considering the complement of a set,
    \[H[\partial E] = -H[\partial (E^c)].\]
    Thanks to this observation, it follows that Lemma~\ref{LEM::sMeanCurvatureLowerBound} also holds true when the sign of the inequalities~\eqref{EQ::MeanCurvatureLowerBound} and~\eqref{EQ::sMeanCurvatureLowerBound} is flipped.
\end{remark}

\begin{figure}[ht]
\centering
\includegraphics[width=0.8\textwidth]{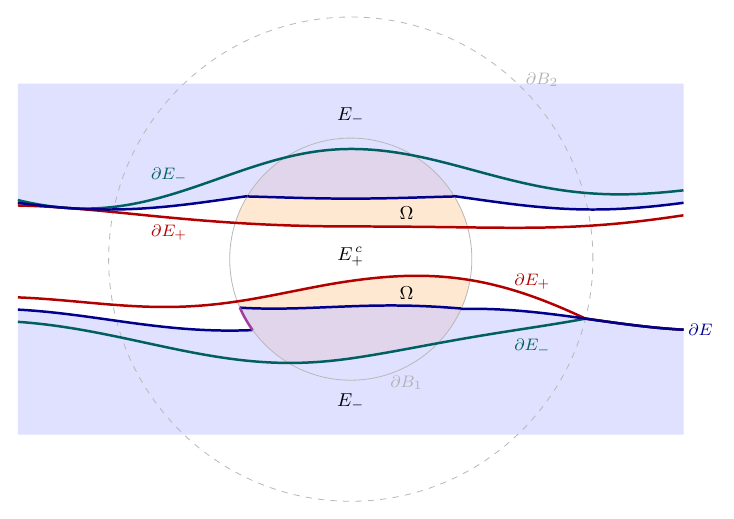}
\caption{Sketch illustrating the situation described by Lemma~\ref{LEM::SlabMinimizerExistence}.}\label{FIG::Slab}
\end{figure}

We now show how the lower bound on the~$s$-mean curvature given by Lemma~\ref{LEM::sMeanCurvatureLowerBound} allows for a Perron-type procedure to take place, yielding a stable local minimizer for the fractional perimeter in the slab contained between two barriers. See also Figure~\ref{FIG::Slab} for a visual representation.

\begin{lemma}\label{LEM::SlabMinimizerExistence}
    Let~$E_+,E_- \subset \R^n$ be such that~$\partial E_\pm \cap B_{2}$ is a~$C^{2,\alpha}$ hypersurface for some~$\alpha \in (0,1)$. Let~$s \in (0,1)$ and suppose that, for every~$x \in \partial E_+ \cap B_1$ and~$y \in \partial E_- \cap B_1$, there holds
    \begin{equation}\label{EQ::sMeanCurvatureBarriers}
        H_s[\partial E_+](x) \geq 1 \qquad \text{and} \qquad H_s[\partial E_-](y) \leq -1.
    \end{equation}
    Furthermore, assume that~$E_- \subset E_+$, denote~$\Omega \coloneq ( E_+ \smallsetminus E_-) \cap B_1$, and let~$E_0 \subset \R^n$ be such that
    \begin{equation}\label{EQ::E0Assumption}
     E_- \subset E_0 \subset E_+ \smallsetminus \Omega.        
    \end{equation}

    Then, there exists~$E \subset \R^n$ such that~$E \smallsetminus \Omega = E_0$, with
    \begin{equation}\label{EQ::SlabMinimality}
        \pers(E, \Omega) = \inf_{\substack{F \subset \R^n\\F \smallsetminus \Omega = E_0}} \pers(F, \Omega)
    \end{equation}
    and
    \begin{equation}\label{EQ::BarrierStrictness}
        \partial E \cap (\partial E_- \cup \partial E_+) \cap B_1 = \varnothing.
    \end{equation}

    Moreover,~$E$ is a stable critical point of~$\pers(\cdot, B_1)$.
\end{lemma}
\begin{proof}
    We claim that
    \begin{equation}\label{EQ::BarriersCannotTouch}
        \partial E_+ \cap \partial E_- \cap B_1 = \varnothing.
    \end{equation}
    Indeed, suppose there exists~$x \in \partial E_+ \cap \partial E_- \cap B_1$. Then, by the inclusion~$E_- \subset E_+$, we have that
    \begin{equation}
        H_s[\partial E_+](x) - H_s[\partial E_-](x) = \pv\int_{\R^n}\frac{ (\chi_{E_+^c}(y) - \chi_{E_-^c}(y))+(\chi_{E_-}(y) - \chi_{E_+}(y))}{|x-y|^{n+s}}\,dy \leq 0.
    \end{equation}
    However, this together with~\eqref{EQ::sMeanCurvatureBarriers}, implies that
    \begin{equation}
  0 < 2 \leq H_s[\partial E_+](x) - H_s[\partial E_-](x) \leq 0,
    \end{equation}
    a contradiction. Thus,~\eqref{EQ::BarriersCannotTouch} is established.

    Since~$\Omega \subset B_1 \Subset B_{\frac{3}{2}}$ and the~$s$-perimeter is non-negative, we know that
    \begin{equation}
        \pers(\Omega) = \pers(\Omega, B_{\frac{3}{2}}) \leq \pers(\Omega, B_{\frac{3}{2}}) + \pers((
        E_+ \smallsetminus E_-) \cup B_1, B_{\frac{3}{2}}),
    \end{equation}
    hence, by submodularity,
    \begin{equation}
        \pers(\Omega) \leq \pers( E_+ \smallsetminus E_-, B_{\frac{3}{2}}) + \pers(B_1, B_{\frac{3}{2}}) < +\infty,
    \end{equation}
    also thanks to the regularity of~$\partial E_\pm \cap B_2$. Therefore, the set~$\Omega$ has finite fractional perimeter, and the existence of a minimizer~$E$ is given by direct minimization (see also~\cite{MR2675483}*{Section~3}).

    Consider now~$x \in \partial E \cap (\partial E_+ \cup \partial E_-) \cap B_1$. Thanks to~\eqref{EQ::BarriersCannotTouch} we have that either~$x \in \partial E_+$ or~$x \in \partial E_-$, but not both. Let us suppose that~$x \in \partial E_-$. Then, we can find~$r > 0$ so small that~$B_r(x) \subset E_+\cap B_1$.
   
   As a result, we find that~$E$ is a variational supersolution
    of the~$s$-perimeter inside~$B_r(x)$, according to
    the definition given in~\cite{MR2675483}*{Section~2} (see in particular formula~(2.1) there). Indeed,
     for every~$A \subset B_r(x) \smallsetminus E$,
     \begin{eqnarray*} &&L(A,E)- L(A, (E\cup A)^c)=
     \pers(E,B_r(x)) - \pers(E \cup A,B_r(x))
     \\&&\qquad=\pers(E,\Omega) - \pers(E \cup A,\Omega)\leq 0,
     \end{eqnarray*}
     thanks to the minimality of~$E$ in~$\Omega$.
    
     Also,~$E$ contains the interior tangent ball to~$E_-$ at~$x$, and therefore, by~\cite{MR2675483}*{Theorem~5.1}, we have the Euler-Lagrange inequality
    \begin{equation}
        H_s[\partial E](x) \geq 0.
    \end{equation}
    The inclusion~$E_- \subset E$ and~\eqref{EQ::sMeanCurvatureBarriers} entail that
    \begin{equation}
 0 > -1 \geq H_s[\partial E_-](x) \geq H_s[\partial E](x) \geq 0,
    \end{equation}
    a clear contradiction, so that~$\partial E \cap \partial E_- \cap B_1$ is empty.

    We handle the case~$x \in \partial E_+$ by symmetry. If we replace~$E_+$,~$ E_-$, and~$E$ with~$E_-^c$,~$E_+^c$, and~$E^c$ respectively, the assumptions of Lemma~\ref{LEM::SlabMinimizerExistence} and~\eqref{EQ::SlabMinimality} remain true. Thus, we see that~$\partial E^c \cap \partial E_+ \cap B_1$ is empty, concluding the proof of~\eqref{EQ::BarrierStrictness}.

    It remains to prove that~$E$ is stable in~$B_1$. Let~$X \in C^\infty_c(B_1)$. Consider the following partition
    \begin{equation}
        \partial \Omega = (\partial E_+ \cap B_1) \cup (\partial E_- \cap B_1) \cup \big(\partial B_1 \cap (E_+ \smallsetminus E_-)\big).
    \end{equation}
    Thanks to this,~\eqref{EQ::BarrierStrictness}, and the fact that~$X$ has compact support, we have that
    \begin{equation}
        \tau \coloneq \operatorname{dist}(\partial E \cap \operatorname{spt}(X), \partial \Omega) > 0.
    \end{equation}
    Hence, by the boundedness of~$|X|$, there exists~$t_\ast > 0$ such that, for every~$t \in (-t_\ast, t_\ast)$, we have
    \begin{equation}
        \bigg(\Phi^X(E, t)\bigg) \triangle E \subset \Omega,
    \end{equation}
    where~$\triangle$ is the symmetric difference of two sets. Namely,~$\Phi^X(E, t)$ is a suitable competitor for the minimality of~$\pers(E, \Omega)$, therefore
    \begin{equation}
        \pers(E, \Omega) \leq \pers(\Phi^X(E, t), \Omega)\qquad\text{for every }t\in(-t_\ast, t_\ast).
    \end{equation}
    Since~$E$ and~$\Phi^X(E, t)$ coincide outside~$\Omega$, we infer that
    \begin{equation}
        \pers(E, B_1) \leq \pers(\Phi^X(E, t), B_1)\qquad\text{for every }t\in(-t_\ast, t_\ast),
    \end{equation}
    from which the desired minimality and stability properties plainly follow.
\end{proof}

With Lemma~\ref{LEM::SlabMinimizerExistence} at hand, we are now in a position to prove Theorem~\ref{THM::Main}.

\subsection{Conclusion of the Proof of Theorem~\ref{THM::Main}}
Let~${\{\Sigma_t\}}_t$ be the foliation given by Theorem~\ref{THM::Foliation} applied to~$\Sigma \coloneq \partial E \cap B_2$. By~\eqref{EQ::eps0OptimalCondition} we can pick
\begin{equation}\label{EQ::tBarDefinition}
    \bar{t} \coloneq \frac{\min\{1,c_1\}\min\{1,\sep(\Sigma)\}}{\max\{1, 2C_\ast\}} < \eps_0,
\end{equation}
where~$c_1$ is~$c$ from~\eqref{EQ::eps0OptimalCondition}, and~$C_\ast$ is as in Lemma~\ref{LEM::FixedPoint}, and consider
\begin{equation}\label{EQ::EplusminusDefinition}
    E_+ \coloneq E \cup \bigcup_{t \in (0,\bar{t})} \Sigma_t \qquad \text{and} \qquad E_- \coloneq E \smallsetminus \bigcup_{t \in (-\bar{t},0)} \Sigma_t.
\end{equation}
Here, $\Sigma$ is oriented by the outer unit normal to~$E$, so that, by the ordering of the foliation, $\Sigma_t$ lies on the side of~$\nu$ for~$t > 0$. As a result,~$\Sigma_t \subset \bar{E}$ if and only if~$t \leq 0$ and, similarly,~$\Sigma_t \subset E^c$ if and only if~$t \geq 0$. Hence, from~\eqref{EQ::EplusminusDefinition} follows that
\begin{equation}
    \partial E_+ \cap B_2 = \Sigma_{\bar{t}} \cap B_2 \qquad \text{and} \qquad \partial E_- \cap B_2 = \Sigma_{-\bar{t}} \cap B_2.
\end{equation}
In particular, this and~\eqref{EQ::FoliationCMCCondition}
entail that, for every~$x \in \partial E_+ \cap B_2$ and~$y \in \partial E_- \cap B_2$, we have that 
\begin{equation}
    H[\partial E_+](x) \equiv \bar{t} \qquad \text{and} \qquad H[\partial E_-](y) \equiv -\bar{t}.
\end{equation}

Applying Lemma~\ref{LEM::sMeanCurvatureLowerBound} with~$\lambda \coloneq \bar{t}$ to~$E_+$, and its variant given by Remark~\ref{REM::sMeanCurvatureLowerBoundSign} to~$E_-$, we find~$s_0 \in (0,1)$ such that, for every~$s \in [s_0,1)$, %chktex 9
\begin{equation}
    H_s[\partial E_+](x) \geq 1 \qquad \text{and} \qquad H_s[\partial E_-](y) \leq -1.
\end{equation}
Moreover, thanks to~\eqref{EQ::optimalS0Estimate} we can pick
\begin{equation}
    s_0 \coloneq  \max\left\{\frac{1}{10},1 - c_2 \min\{1,\rho(\partial E \cap B_2)\} \bar{t}\right\},
\end{equation}
where~$c_2$ is~$c$ from~\eqref{EQ::optimalS0Estimate}. Recalling~\eqref{EQ::tBarDefinition}, we see that
\begin{equation}
    \begin{split}
        s_0 &\leq 1 - c_2 \min\{1,\rho_{\loc}(\partial E \cap B_2), \sep(\Sigma)\} \frac{\min\{1,c_1\}}{\max\{1, 2C_\ast\}} \min\{1,\sep(\Sigma)\}\\
            &\leq 1 - c_2 \min\{1,\rho_{\loc}(\partial E \cap B_2)\} \frac{\min\{1,c_1\}}{\max\{1, 2C_\ast\}} (\min\{1,\sep(\Sigma)\})^2\\
            &= 1 - c_2 \min\{1,\rho_{\loc}(\partial E \cap B_2)\} \frac{\min\{1,c_1\}}{\max\{1, 2C_\ast\}} \min\{1,\sep^2(\Sigma)\}.
    \end{split}
\end{equation}
Hence, for every~$s \in (0,1)$ such that~\eqref{EQ::s0Constraint} holds true with
\[C \coloneq c_2 \min\{1,\rho_{\loc}(\partial E \cap B_2)\} \frac{\min\{1,c_1\}}{\max\{1, 2C_\ast\}},\]
we also have that~$s \geq s_0$. As a consequence, for every such~$s$, we can apply Lemma~\ref{LEM::SlabMinimizerExistence} with~$E_0 \coloneq E \smallsetminus \Omega$ to find~$F \subset \R^n$ such that~$\partial F$ is a stable nonlocal minimal surface in~$B_1$.

\begin{figure}[ht]
\centering
\includegraphics[width=0.5\textwidth]{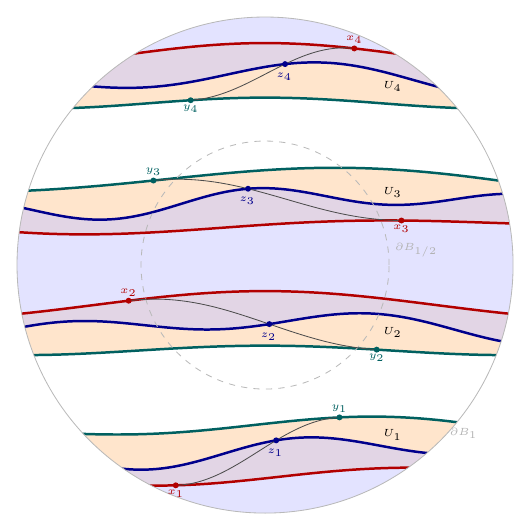}
\caption{Sketch of the points~$x_j$,~$y_j$, and~$z_j$, with~$j \in \{1,2,3,4\}$, in the proof of Theorem~\ref{THM::Main}.}\label{FIG::kComponents}
\end{figure}

It remains to prove that~$\partial F \cap B_1$ has at least~$k$ connected components. Let~${\{\Sigma^j\}}_{j=1}^k$ denote the distinct connected components of~$\partial E \cap B_{\frac{1}{2}}$. Let
\begin{equation}
    U_j \coloneq \bigcup_{t \in (-\bar{t},\bar{t})} \Sigma_t^j.
\end{equation}
Clearly, since~$\bar{t} \leq \sep(\Sigma)$, we have that
\begin{equation}\label{EQ::UjDisconnected}
    \operatorname{dist}\left(U_j, U_i\right) > 0
\end{equation}
whenever~$1\leq i < j \leq k$. Moreover, since~$\bar{t} < \frac{1}{C_\ast}$, for every~$j \in \{1,\dots,k\}$, the set~$\overline{U_j}$ is contained in~$B_1$. Namely, there exist~$x_j \in \Sigma_{ -\bar{t}}^j \cap B_1$ and~$y_j \in \Sigma_{\bar{t}}^j \cap B_1$. Importantly,~$x_j \in \partial E_-$ and~$y_j \in \partial E_+$. Therefore,~$x_j \in F$ and~$y_j \notin F$. Exploiting the path connectedness of~$\overline{U_j}$, we find~$z_j \in \partial F \cap U_j$ (see Figure~\ref{FIG::kComponents}). Thanks to~\eqref{EQ::UjDisconnected}, we see that
\[\partial F \cap B_1 \subset (E_+ \smallsetminus E_-) \cap B_1 = \bigcup_{t \in (-\bar{t},\bar{t})} \Sigma_t = \bigcup_{t \in (-\bar{t},\bar{t})} \bigcup_{j=1}^k \Sigma_t^j = \bigcup_{j=1}^k U_j.\]
Therefore, each point~$z_j$ belongs to a different connected component of~$\partial F \cap B_1$, as desired.
\qed{}
\begin{appendix}
\section{CMC foliation around a flat disk}\label{APP::PlaneExample}
In this appendix, we make some considerations regarding Theorem~\ref{THM::Foliation} in the special case where~$\Sigma$ is the flat disk. This provides a concrete example of our result, also allowing us to partially discuss the optimality of~\eqref{EQ::eps0OptimalCondition}.

Consider the flat disk
\begin{equation}
    \Sigma \coloneq \left\{x = (x^\prime,x_n) \in \R^{n-1} \times \R \colon |x^\prime| \leq 2,\, x_n = 0\right\}.
\end{equation}
Clearly,~$\Sigma$ is a minimal surface in~$B_2$ since it has null mean curvature everywhere. Moreover,~$\mathcal{L}_\Sigma$ is just the Laplacian in~$\R^{n-1}$. By identifying~$\Sigma$ with~$B_2 \subset \R^{n-1}$, it is straightforward to see that~$\Sigma$ is strictly stable thanks to the Poincaré inequality.

Let~$R \in (2,+\infty)$ and consider the ball~$B \coloneq B_R(\hat{x})$, where~$\hat{x} \coloneq (0,-\sqrt{R^2-4}) \in \R^{n-1} \times \R$. Clearly, the mean curvature on~$\partial B$ with respect to the outer normal is identically equal to~${\frac{n-1}{R}}$. In particular, consider the function~${u_{R}} \colon \Sigma \to \R$ defined by
\begin{equation}\label{EQ::uRDefinition}
    {u_{R}} \colon x^\prime \longmapsto \sqrt{R^2-|x^\prime|^2} - \sqrt{R^2-4}.
\end{equation}
Using the notation of~\eqref{EQ::SigmaUDefinition}, one can see that~$\Sigma_{{u_{R}}} = \partial B \cap (\R^{n-1} \times \R^+)$. Therefore, the normal geodesic graph of~$u$ over~$\Sigma$ has constant mean curvature equal to~${\frac{n-1}{R}}$. Similarly, exploiting symmetry with respect to the reflection~$x_n \mapsto -x_n$, we have that the normal geodesic graph of~${u_{-R} \coloneq -u_{R}}$ over~$\Sigma$, which is contained in~$\R^{n-1} \times \R^-$, has constant mean curvature equal to~${-\frac{n-1}{R}}$.

\begin{figure}[ht]
\centering
\includegraphics[width=0.5\textwidth]{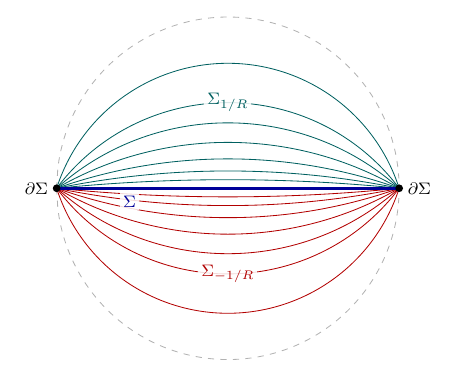}
\caption{The foliation around a flat disk, as provided by Theorem~\ref{THM::Foliation}.}\label{FIG::plane}
\end{figure}

We note that the family~${\{\Sigma_t\}}_{t \in \left({-\frac{n-1}{2},\frac{n-1}{2}}\right)}$, where~$\Sigma_t \coloneq \Sigma_{u_{{\frac{n-1}{t}}}}$, is the foliation around~$\Sigma$ constructed by applying Theorem~\ref{THM::Foliation} and choosing~$\eps_0 \coloneq {\frac{n-1}{2}}$, see Figure~\ref{FIG::plane}. Indeed, from~\eqref{EQ::uRDefinition} it also follows that the leaves enjoy ordering with respect to~$t = {\frac{n-1}{R}}$. In fact, it is well known that there exists no function on~$B_2 \subset \R^{n-1}$ whose graph has constant mean curvature larger than~${\frac{n-1}{2}}$ (the limit case would be the upper and lower halves of~$\partial B_2$), so~$\eps_0$ has a finite optimal value even if~$\Sigma$ is flat and connected, hence~$\sep(\Sigma) = +\infty$ (actually, in this flat case it also has infinite range). This consideration offers an interpretation of the form of~\eqref{EQ::eps0OptimalCondition} and why it cannot be improved to ``$\eps_0 \geq c\sep(\Sigma)$''.

To further understand the implicit limitations that restrict the attainable values of mean curvature for leaves of the foliation, one needs to look into the boundary behavior of the family of functions~$u_{\pm R}$. Since each one of them vanishes on~$\partial \Sigma$, but they are ordered inside, it is interesting to look at their normal (relative to the~$\R^{n-1}$ topology) derivative on~$\partial \Sigma$. Note that, at each~$x^\prime \in \partial \Sigma$, the relative outer normal~$n_\Sigma(x^\prime)$ is equal to~$\frac{x^\prime}{2}$. A simple computation shows that
\begin{equation}
    \partial_{n_\Sigma} {u_{R}}(x^\prime) \equiv -\frac{2}{\sqrt{R^2-4}}.
\end{equation}
In particular, as~$R \searrow 2$, we have that~$|\partial_{n_\Sigma} {u_{R}}(x^\prime)| \to \infty$, showing that the construction cannot be extended further.
\end{appendix}
{
 \hbadness=10000
 \vbadness=10000
 \bibliography{biblio}

@article {MR3107529,
    AUTHOR = {Caffarelli, Luis and Valdinoci, Enrico},
     TITLE = {Regularity properties of nonlocal minimal surfaces via
              limiting arguments},
   JOURNAL = {Adv. Math.},
  FJOURNAL = {Advances in Mathematics},
    VOLUME = {248},
      YEAR = {2013},
     PAGES = {843--871},
      ISSN = {0001-8708,1090-2082},
   MRCLASS = {49Q05 (58E20)},
  MRNUMBER = {3107529},
MRREVIEWER = {Giandomenico\ Orlandi},
       DOI = {10.1016/j.aim.2013.08.007},
       URL = {https://doi.org/10.1016/j.aim.2013.08.007},
}

@article {MR2675483,
    AUTHOR = {Caffarelli, L. and Roquejoffre, J.-M. and Savin, O.},
     TITLE = {Nonlocal minimal surfaces},
   JOURNAL = {Comm. Pure Appl. Math.},
  FJOURNAL = {Communications on Pure and Applied Mathematics},
    VOLUME = {63},
      YEAR = {2010},
    NUMBER = {9},
     PAGES = {1111--1144},
      ISSN = {0010-3640,1097-0312},
   MRCLASS = {49Q05 (58E12)},
  MRNUMBER = {2675483},
MRREVIEWER = {Silvano\ Delladio},
       DOI = {10.1002/cpa.20331},
       URL = {https://doi.org/10.1002/cpa.20331},
}

@article {MR62214,
    AUTHOR = {De Giorgi, Ennio},
     TITLE = {Su una teoria generale della misura {$(r-1)$}-dimensionale in
              uno spazio ad {$r$} dimensioni},
   JOURNAL = {Ann. Mat. Pura Appl. (4)},
  FJOURNAL = {Annali di Matematica Pura ed Applicata. Serie Quarta},
    VOLUME = {36},
      YEAR = {1954},
     PAGES = {191--213},
      ISSN = {0003-4622},
   MRCLASS = {27.2X},
  MRNUMBER = {62214},
MRREVIEWER = {L.\ C.\ Young},
       DOI = {10.1007/BF02412838},
       URL = {https://doi.org/10.1007/BF02412838},
}

@article {MR74499,
    AUTHOR = {De Giorgi, Ennio},
     TITLE = {Nuovi teoremi relativi alle misure {$(r-1)$}-dimensionali in
              uno spazio ad {$r$} dimensioni},
   JOURNAL = {Ricerche Mat.},
  FJOURNAL = {Ricerche di Matematica},
    VOLUME = {4},
      YEAR = {1955},
     PAGES = {95--113},
      ISSN = {0035-5038},
   MRCLASS = {28.0X},
  MRNUMBER = {74499},
MRREVIEWER = {L.\ C.\ Young},
}

@book {MR775682,
    AUTHOR = {Giusti, Enrico},
     TITLE = {Minimal surfaces and functions of bounded variation},
    SERIES = {Monographs in Mathematics},
    VOLUME = {80},
 PUBLISHER = {Birkh\"auser Verlag, Basel},
      YEAR = {1984},
     PAGES = {xii+240},
      ISBN = {0-8176-3153-4},
   MRCLASS = {58E12 (49F10 53A10)},
  MRNUMBER = {775682},
MRREVIEWER = {Helmut\ Kaul},
       DOI = {10.1007/978-1-4684-9486-0},
       URL = {https://doi.org/10.1007/978-1-4684-9486-0},
}

@book {MR2780140,
    AUTHOR = {Colding, Tobias Holck and Minicozzi, II, William P.},
     TITLE = {A course in minimal surfaces},
    SERIES = {Graduate Studies in Mathematics},
    VOLUME = {121},
 PUBLISHER = {American Mathematical Society, Providence, RI},
      YEAR = {2011},
     PAGES = {xii+313},
      ISBN = {978-0-8218-5323-8},
   MRCLASS = {53A10 (35J93 49Q05)},
  MRNUMBER = {2780140},
MRREVIEWER = {Andrew\ Bucki},
       DOI = {10.1090/gsm/121},
       URL = {https://doi.org/10.1090/gsm/121},
}

@article {MR562550,
    AUTHOR = {Fischer-Colbrie, Doris and Schoen, Richard},
     TITLE = {The structure of complete stable minimal surfaces in
              {$3$}-manifolds of nonnegative scalar curvature},
   JOURNAL = {Comm. Pure Appl. Math.},
  FJOURNAL = {Communications on Pure and Applied Mathematics},
    VOLUME = {33},
      YEAR = {1980},
    NUMBER = {2},
     PAGES = {199--211},
      ISSN = {0010-3640,1097-0312},
   MRCLASS = {53C40 (58E12)},
  MRNUMBER = {562550},
MRREVIEWER = {Themistocles\ M.\ Rassias},
       DOI = {10.1002/cpa.3160330206},
       URL = {https://doi.org/10.1002/cpa.3160330206},
}

@article {MR2782803,
    AUTHOR = {Caffarelli, Luis and Valdinoci, Enrico},
     TITLE = {Uniform estimates and limiting arguments for nonlocal minimal
              surfaces},
   JOURNAL = {Calc. Var. Partial Differential Equations},
  FJOURNAL = {Calculus of Variations and Partial Differential Equations},
    VOLUME = {41},
      YEAR = {2011},
    NUMBER = {1-2},
     PAGES = {203--240},
      ISSN = {0944-2669,1432-0835},
   MRCLASS = {49Q05 (28A75 58E12)},
  MRNUMBER = {2782803},
MRREVIEWER = {Andrew\ Bucki},
       DOI = {10.1007/s00526-010-0359-6},
       URL = {https://doi.org/10.1007/s00526-010-0359-6},
}

@article {MR2765717,
    AUTHOR = {Ambrosio, Luigi and De Philippis, Guido and Martinazzi, Luca},
     TITLE = {Gamma-convergence of nonlocal perimeter functionals},
   JOURNAL = {Manuscripta Math.},
  FJOURNAL = {Manuscripta Mathematica},
    VOLUME = {134},
      YEAR = {2011},
    NUMBER = {3-4},
     PAGES = {377--403},
      ISSN = {0025-2611,1432-1785},
   MRCLASS = {49Q15 (49J45 49Q20)},
  MRNUMBER = {2765717},
MRREVIEWER = {Matteo\ Focardi},
       DOI = {10.1007/s00229-010-0399-4},
       URL = {https://doi.org/10.1007/s00229-010-0399-4},
}

@article {MR4581189,
    AUTHOR = {Dipierro, Serena and Valdinoci, Enrico},
     TITLE = {Some perspectives on (non)local phase transitions and minimal
              surfaces},
   JOURNAL = {Bull. Math. Sci.},
  FJOURNAL = {Bulletin of Mathematical Sciences},
    VOLUME = {13},
      YEAR = {2023},
    NUMBER = {1},
     PAGES = {Paper No. 2330001, 77},
      ISSN = {1664-3607,1664-3615},
   MRCLASS = {53A10 (35J61 35R11 49Q05 82B26)},
  MRNUMBER = {4581189},
MRREVIEWER = {Nicola\ Abatangelo},
       DOI = {10.1142/S1664360723300013},
       URL = {https://doi.org/10.1142/S1664360723300013},
}

@article {MR3230079,
    AUTHOR = {Abatangelo, Nicola and Valdinoci, Enrico},
     TITLE = {A notion of nonlocal curvature},
   JOURNAL = {Numer. Funct. Anal. Optim.},
  FJOURNAL = {Numerical Functional Analysis and Optimization. An
              International Journal},
    VOLUME = {35},
      YEAR = {2014},
    NUMBER = {7-9},
     PAGES = {793--815},
      ISSN = {0163-0563,1532-2467},
   MRCLASS = {53A99},
  MRNUMBER = {3230079},
MRREVIEWER = {Zhenhai\ Liu},
       DOI = {10.1080/01630563.2014.901837},
       URL = {https://doi.org/10.1080/01630563.2014.901837},
}

@article {MR110078,
    AUTHOR = {Federer, Herbert},
     TITLE = {Curvature measures},
   JOURNAL = {Trans. Amer. Math. Soc.},
  FJOURNAL = {Transactions of the American Mathematical Society},
    VOLUME = {93},
      YEAR = {1959},
     PAGES = {418--491},
      ISSN = {0002-9947,1088-6850},
   MRCLASS = {53.00},
  MRNUMBER = {110078},
MRREVIEWER = {C.\ B.\ Allendoerfer},
       DOI = {10.2307/1993504},
       URL = {https://doi.org/10.2307/1993504},
}

@article {MR2481045,
    AUTHOR = {Pacard, F. and Xu, X.},
     TITLE = {Constant mean curvature spheres in {R}iemannian manifolds},
   JOURNAL = {Manuscripta Math.},
  FJOURNAL = {Manuscripta Mathematica},
    VOLUME = {128},
      YEAR = {2009},
    NUMBER = {3},
     PAGES = {275--295},
      ISSN = {0025-2611,1432-1785},
   MRCLASS = {53C42},
  MRNUMBER = {2481045},
MRREVIEWER = {Alma\ Luisa\ Albujer Brotons},
       DOI = {10.1007/s00229-008-0230-7},
       URL = {https://doi.org/10.1007/s00229-008-0230-7},
}

@book {MR2759829,
    AUTHOR = {Brezis, Haim},
     TITLE = {Functional analysis, {S}obolev spaces and partial differential
              equations},
    SERIES = {Universitext},
 PUBLISHER = {Springer, New York},
      YEAR = {2011},
     PAGES = {xiv+599},
      ISBN = {978-0-387-70913-0},
   MRCLASS = {35-01 (46-01 46E35 46N20 47F05)},
  MRNUMBER = {2759829},
MRREVIEWER = {Vicen\c tiu\ D.\ R\u adulescu},
}

@article{chan2025nonlocalapproximationminimalsurfaces,
      title={Nonlocal approximation of minimal surfaces: optimal estimates from stability}, 
      author={Hardy Chan and Serena Dipierro and Joaquim Serra and Enrico Valdinoci},
       JOURNAL = {Comm. Pure Appl. Math.},
      year={2025},
      eprint={2308.06328},
      archivePrefix={arXiv},
      primaryClass={math.DG},
      url={https://arxiv.org/abs/2308.06328}, 
}

@article {MR3798717,
    AUTHOR = {D\'avila, Juan and del Pino, Manuel and Wei, Juncheng},
     TITLE = {Nonlocal {$s$}-minimal surfaces and {L}awson cones},
   JOURNAL = {J. Differential Geom.},
  FJOURNAL = {Journal of Differential Geometry},
    VOLUME = {109},
      YEAR = {2018},
    NUMBER = {1},
     PAGES = {111--175},
      ISSN = {0022-040X,1945-743X},
   MRCLASS = {53A10 (49Q05 58E12)},
  MRNUMBER = {3798717},
MRREVIEWER = {Chun-Chi\ Lin},
       DOI = {10.4310/jdg/1525399218},
       URL = {https://doi.org/10.4310/jdg/1525399218},
}

@article {MR1084717,
    AUTHOR = {Ye, Rugang},
     TITLE = {Foliation by constant mean curvature spheres},
   JOURNAL = {Pacific J. Math.},
  FJOURNAL = {Pacific Journal of Mathematics},
    VOLUME = {147},
      YEAR = {1991},
    NUMBER = {2},
     PAGES = {381--396},
      ISSN = {0030-8730,1945-5844},
   MRCLASS = {53C12},
  MRNUMBER = {1084717},
MRREVIEWER = {Alexandre\ L.\ Yampol\cprime ski\u i},
       URL = {http://projecteuclid.org/euclid.pjm/1102644918},
}

@article {MR1101226,
    AUTHOR = {White, Brian},
     TITLE = {The space of minimal submanifolds for varying {R}iemannian
              metrics},
   JOURNAL = {Indiana Univ. Math. J.},
  FJOURNAL = {Indiana University Mathematics Journal},
    VOLUME = {40},
      YEAR = {1991},
    NUMBER = {1},
     PAGES = {161--200},
      ISSN = {0022-2518,1943-5258},
   MRCLASS = {58D10 (53C42)},
  MRNUMBER = {1101226},
MRREVIEWER = {Jo\~ao\ Lucas Marques Barbosa},
       DOI = {10.1512/iumj.1991.40.40008},
       URL = {https://doi.org/10.1512/iumj.1991.40.40008},
}

@article {MR3516886,
    AUTHOR = {Dipierro, Serena and Savin, Ovidiu and Valdinoci, Enrico},
     TITLE = {Graph properties for nonlocal minimal surfaces},
   JOURNAL = {Calc. Var. Partial Differential Equations},
  FJOURNAL = {Calculus of Variations and Partial Differential Equations},
    VOLUME = {55},
      YEAR = {2016},
    NUMBER = {4},
     PAGES = {Art. 86, 25},
      ISSN = {0944-2669,1432-0835},
   MRCLASS = {49Q05 (35R11 53A10)},
  MRNUMBER = {3516886},
MRREVIEWER = {Martin\ Man-Chun\ Li},
       DOI = {10.1007/s00526-016-1020-9},
       URL = {https://doi.org/10.1007/s00526-016-1020-9},
}

@incollection {MR3824212,
    AUTHOR = {Dipierro, Serena and Valdinoci, Enrico},
     TITLE = {Nonlocal minimal surfaces: interior regularity, quantitative
              estimates and boundary stickiness},
 BOOKTITLE = {Recent developments in nonlocal theory},
     PAGES = {165--209},
 PUBLISHER = {De Gruyter, Berlin},
      YEAR = {2018},
      ISBN = {978-3-11-057156-1; 978-3-11-057155-4},
   MRCLASS = {49Q05 (35A15)},
  MRNUMBER = {3824212},
MRREVIEWER = {Annalisa\ Cesaroni},
       DOI = {10.1515/9783110571561-006},
       URL = {https://doi.org/10.1515/9783110571561-006},
}

@article {MR3981295,
    AUTHOR = {Cinti, Eleonora and Serra, Joaquim and Valdinoci, Enrico},
     TITLE = {Quantitative flatness results and {$BV$}-estimates for stable
              nonlocal minimal surfaces},
   JOURNAL = {J. Differential Geom.},
  FJOURNAL = {Journal of Differential Geometry},
    VOLUME = {112},
      YEAR = {2019},
    NUMBER = {3},
     PAGES = {447--504},
      ISSN = {0022-040X,1945-743X},
   MRCLASS = {49Q05 (35R11 53A10)},
  MRNUMBER = {3981295},
MRREVIEWER = {Doan The Hieu},
       DOI = {10.4310/jdg/1563242471},
       URL = {https://doi.org/10.4310/jdg/1563242471},
}

@article {MR4116635,
    AUTHOR = {Cabr\'e, Xavier and Cinti, Eleonora and Serra, Joaquim},
     TITLE = {Stable {$s$}-minimal cones in {$\mathbb{R}^3$} are flat for
              {$s\sim 1$}},
   JOURNAL = {J. Reine Angew. Math.},
  FJOURNAL = {Journal f\"ur die Reine und Angewandte Mathematik. [Crelle's
              Journal]},
    VOLUME = {764},
      YEAR = {2020},
     PAGES = {157--180},
      ISSN = {0075-4102,1435-5345},
   MRCLASS = {53A10},
  MRNUMBER = {4116635},
MRREVIEWER = {Yongsheng\ Zhang},
       DOI = {10.1515/crelle-2019-0005},
       URL = {https://doi.org/10.1515/crelle-2019-0005},
}
}
\end{document}